\documentclass[twocolumn]{autart}

\usepackage{natbib}
\usepackage{mathrsfs}
\usepackage{eufrak}
\setcitestyle{authoryear,open={(},close={)}}
\usepackage{mathtools}
\usepackage{grffile}
\usepackage{cases}
\usepackage{bm}
\usepackage[english]{babel}
\usepackage[autostyle]{csquotes}
\usepackage{multirow}
\usepackage{caption}
\usepackage{graphicx}
\usepackage{wrapfig}
\usepackage{amssymb,amsmath}
\usepackage{epsfig}
\usepackage{epstopdf}
\usepackage{setspace}
\usepackage{enumerate}
\usepackage{tikz}
\usepackage{mdwlist}
\usepackage{color}
\usepackage{algorithm}
\usepackage{algcompatible}
\usepackage{float}
\usepackage{algpseudocode}
\usepackage{pifont}
\usepackage{url}
\usepackage[makeroom]{cancel}
\usetikzlibrary{shapes,arrows}
\allowdisplaybreaks
\usepackage[normalem]{ulem} 

\newtheorem{example}{Example}
\newtheorem{theorem}{Theorem}
\newtheorem{proposition}{Proposition}
\newtheorem{lemma}{Lemma}
\newtheorem{remark}{Remark}

\newtheorem{assumption}{Assumption}

\newtheorem{definition}{Definition}

\newcommand{\D}{{\Delta}}
\newcommand{\R}{{\mathbb{R}}}

\newcommand{\bx}{{\mathbf{x}}}

\newcommand{\bq}{{\mathbf{q}}}

\newcommand{\bz}{{\mathbf{z}}}

\newcommand{\one}{{\mathbf{1}}}
\newcommand{\zero}{{\mathbf{0}}}
\newcommand{\CG}{{\Phi}}

\begin{document}
\begin{frontmatter}

\title{Distributed Nash Equilibrium Seeking under Partial-Decision Information via the Alternating Direction Method of Multipliers  
} 

\vspace{-0.4cm}
\author[First]{Farzad Salehisadaghiani},
\author[Second]{Wei Shi}, 
\author[First]{\hspace{0.2cm}Lacra Pavel}

\address[First]{Department of Electrical and Computer Engineering, University of Toronto, Toronto, ON M5S 3G4, Canada \\{\tt\small farzad.salehisadaghiani@mail.utoronto.ca, pavel@control.utoronto.ca}.}
\address[Second]{School of Electrical, Computer and Energy Engineering, Arizona State University, Tempe, AZ 85287,  USA \\{\tt\small  wilbur.shi@asu.edu}.} 

\begin{abstract}       
In this paper, we consider the problem of finding a Nash equilibrium in a multi-player game over generally connected networks. This model differs from a conventional setting in that players have partial information on the actions of their opponents and  the communication graph is not necessarily the same as the players' cost dependency graph. We develop a relatively fast algorithm within the framework of inexact-ADMM, based on local  information exchange between the players. We prove  its convergence to Nash equilibrium  for fixed step-sizes and analyze its  convergence rate. Numerical simulations illustrate its benefits when compared to a consensus-based gradient type algorithm with diminishing step-sizes.
\end{abstract}



\end{frontmatter}

\section{Introduction}
\vspace{-0.4cm}
 
{ 
We consider distributed Nash equilibrium (NE)  seeking in a setting where players  have limited local information, over a communication network. In contrast to  the classical setting where all players have access to all their opponents' decisions  \cite{Scutari2014}, \cite{frihauf2012nash}, \cite{shanbhag}, our interest lies in networked regimes where agents may only access or observe the decisions of their neighbours, and   
 there is no central node that has bidirectional communication with all players. 
{  This is a topic of increasing  recent interest given the proliferation of engineering networked applications requiring distributed protocols that operate under partial, local information (V2V, P2P, ad-hoc, smart-grid networks), and the deregulation of global markets.}  Applications that motivate us to consider NE seeking in such a setting range from 
spectrum access in cognitive radio networks, \cite{WangWu2010}, \cite{Cheng2014}, \cite{swenson2015empirical},  rate control and congestion games in wireless, vehicular or ad-hoc  networks, \cite{Alpcan2005}, \cite{yin2011nash}, \cite{Tekin2012}, to networked Nash-Cournot competition, \cite{Bimpikis2014}, \cite{Jayash}, and opinion dynamics in social networks, \cite{Srikant2014}, \cite{Ozdaglar2016}. 
For example, in cognitive radio networks users adaptively adjust their operating parameters based on interactions with the environment and other users in the network,  \cite{WangWu2010}. In a vehicular ad-hoc network setting, vehicles on the road often maintain a relatively stable topology and 
form clusters wherein vehicles can communicate one with another. 
A best-response algorithm  for channel selection is proposed in \cite{Cheng2014}, based on vehicles within a cluster communicating with each other, while selfishly trying to maximize their own utilities. A similar information exchange setting is considered in \cite{swenson2015empirical}, where because players  are unable to directly observe the actions of all others, they  engage in local, non-strategic information exchange to eventually agree on a common Nash equilibrium. 
{In an economic setting,  \cite{bramoulle2014strategic}
 draws attention to the problem of ``who interacts with whom" in a network and to the importance of communication with neighbouring players. An example  is Nash-Cournot competition played over a network,  between 
a set of firms that compete over a set of locations, \cite{Bimpikis2014}. A communication network is formed between the firms, and prescribes how they communicate locally their production decision over those locations,  \cite{Jayash}. 
All these examples are  non-cooperative in the way decisions are made (each agent is self-interested in minimizing only its own cost), but collaborative in information sharing (agents have the incentive to 
exchange information with their neighbours in order to  mitigate the lack of global information on others' decisions). 
Unlike the classical setting, where all others' decisions are known and information exchange is not incentive compatible, in a partial-decision information setting (not usually considered in classical game theory), information exchange  is motivated because of the limited data each player has. 
Such an information exchange setting has started to be considered very recently,  \cite{swenson2015empirical},  \cite{Jayash}, \cite{salehisadaghiani2016distributed}, but convergence is guaranteed only for diminishing step-sizes.  
\textit{Motivated by the above, in this paper we develop a distributed NE seeking algorithm that operates under limited (partial) decision-information, and is guaranteed to converge to NE with constant step-sizes. Each player updates his decision and his estimate simultaneously in a single step, based on local information exchange with his neighbours.} } 
 }  
 \vspace{-0.2cm}
 
 We use an  ADMM approach to develop the algorithm.  Originally  developed in 1970s,  the \emph{Alternating Direction Method of Multipliers} (ADMM) method has become widely used in distributed optimization problems (DOP) after its re-introduction in \cite{boyd2011distributed}, \cite{wei2012distributed}, \cite{WeiOzdaglar2013arX}, \cite{shi2014linear}, \cite{chang2015multi}, \cite{HongLuoRaza2016}. 
In a DOP, $N$ agents  cooperatively minimize a global objective, 
$f(y):=\sum_{i=1}^Nf_i(y)$,  over $y \in \Omega$, by communicating over a connected graph. 
The problem is reformulated by introducing a separate decision variable $x_i$ for each agent $i$ and imposing the equality constraint $x_i=x_j$ for all agents connected by an edge (neighbours),  leading to  \vspace{-0.3cm} 
\begin{equation}\label{miniBis}
\begin{cases}
\begin{aligned}
& \underset{x}{\text{minimize}}
& & f(x):=\sum_{i=1}^N f_i(x_i) \\
& \text{subject to}
& & A x =0, \quad x_i \in \Omega \, \, \, \forall  i =1, \ldots,N
\end{aligned}
\end{cases}
\end{equation} 
where   $x=[x_i]^N_{i=1}$, $A$ is the edge-node incidence matrix. 
{ADMM for DOP relies on the additive structure and separability of the objective function \eqref{miniBis} and on the linearity of the equality constraints. Each agent $i$ will minimize its cost $f_i(x_i)$ over  $x_i$ and this $f_i$ is independent (\emph{decoupled}) of the other agents'  $x_j$, $j \neq i$. 
}Dual decomposition leads to $N$ parallel decoupled dual ascent problems. 
\vspace{-0.255cm}

In this work, we exploit the benefits of ADMM in the context of finding a NE of a game, where 
each player (agent) $i$ aims to minimize its own cost function $J_i$ with respect to (w.r.t) its action $x_i$,  given a profile of other players' actions,  except himself, $x_{-i}$, 
\vspace{-0.3cm}
\begin{equation}
\label{mini_00}
\begin{cases}
\begin{aligned}
& \underset{x_i}{\text{minimize}}
& & J_i(x_i,x_{-i}) \\
& \text{subject to}
& & x_i\in \Omega_i
\end{aligned}
\end{cases}\quad\forall i =1, \ldots,N.
\end{equation}
We consider  a networked information exchange setting as in the motivating examples, due to players' limited (partial) information on the others' decisions. 
There are several  challenges when comparing a game to a DOP:\vspace{-0.255cm}
\begin{itemize}
	\item A Nash game can be seen as a set of parallel \emph{coupled} optimization problems,  \eqref{mini_00}. Each player's cost is dependent on the other players' decisions, hence its decision is directly affected by these. 
	\item Each player $i$ updates only its own decision $x_i$, but he  also requires an estimate of all others $[x_j]_{j =1,\ldots,N,\,j\neq i}$,  or $x_{-i}$, in order to solve its optimization problem. 
\end{itemize} \vspace{-0.25cm}
We introduce local estimates of players' actions and (virtual) constraints for their consensus. We reformulate the problem and relate it to a modified game (a set of optimization problems) with  consensus constraints.  
As in classical ADMM we reduce the computational costs by using a linear approximation in players' action update rule (inexact-ADMM). {However, direct application of ADMM is not possible because each minimization is not over the whole vector (estimate), but rather over part of it (action).}  Every player updates its \emph{decision/action} as well as its \emph{estimates} of the other players' actions by this synchronous, inexact ADMM-type algorithm, using its own estimates and those of its neighbours. We show that the algorithm converges to  NE under fixed step-sizes, the first such algorithm to the best of our knowledge. 
{A short version 
appears in \cite{salehiPavelIFAC2017}, \cite{WeiShiACC} under stronger assumptions. }

\vspace{-0.255cm}

\emph{\textbf{Related Works.}} 
Distributed NE computation  is relevant to many  applications, such as sensor network coordination \cite{Stankovic2012}, flow control \cite{Alpcan2005}, optical networks, \cite{jayash11}. {In the \emph{classical setting} of distributed NE computation  using best-response or gradient-based schemes,  competitors' decisions are assumed to be observable by all players, \cite{Facchi24}, \cite{frihauf2012nash}. }
Recent work  in the classical setting considered games with monotone pseudo-gradient \cite{Scutari2014,shanbhag,Yousefian2015}, \cite{zhu2016distributed},  networked aggregative games with quadratic cost functions, \cite{parise2015network}, or games with affine constraints,  \cite{yipengLP1, peng2}. Players are assumed to have access to the others' decisions, so the issue of \emph{partial information on the opponents' decisions} is not considered in the above works. 

\vspace{-0.255cm} However, having access to all others' decisions can be impractical in distributed  networks, \cite{Marden2007_2}, \cite{Marden2013}. 
{In recent years, there has been an increasing interest  to 
consider this issue and how to deal with it via networked (partial)-information exchange. 
Our work is related to this literature. }
Consensus-based approaches have been proposed,  by using  fictitious play in congestion games,  \cite{swenson2015empirical}, 
or projected-gradient algorithms, either for the special class of  \emph{aggregative} games, \cite{Jayash}, or for games where players' cost functions depend on others' actions in a general manner, in \cite{salehisadaghiani2016distributed}, or are partially coupled, \cite{salehisadaghianiPavel2018Auto}. Convergence to the NE was shown only for diminishing step-sizes, under strict monotonicity of the pseudo-gradient. 
In contrast, our ADMM algorithm converges with fixed step-sizes and has faster convergence. 
We note that very recently continuous-time NE seeking dynamics are proposed based on consensus and gradient-type dynamics, for games with unconstrained action sets, \cite{YeHu2017}, or based on passivity, \cite{GadjovPavelTAC2018}. 
 Differently,  we develop a discrete-time NE seeking algorithm  for games with compact action sets, based on an ADMM approach.

\vspace{-0.255cm}
{Our work is related to the large literature on ADMM for DOP but there are several differences. In contrast to typical ADMM algorithms developed for DOP \eqref{miniBis} which is separable, e.g.  \cite{wei2012distributed}, \cite{WeiOzdaglar2013arX}, \cite{shi2014linear}, \cite{chang2015multi},  \cite{HongLuoRaza2016}, our ADMM algorithm for NE seeking of \eqref{mini_00} has an extra step for estimate update, besides the update of each player's decision variable. This is because of the intrinsic coupling in $J_i$ \eqref{mini_00}; the update of each player's action is  coupled nonlinearly to this estimate and to the others' decisions. This is unlike ADMM for DOP where, due to decoupled costs, the update of each agent's decision is only linearly coupled to the others' decisions (via the constraints). This also leads to a technical difference on the convexity assumptions under which convergence is shown.
Typical assumptions in DOP are individual (strict/strong) joint convexity of each agent's decoupled cost function $f_i$  \eqref{miniBis},  \cite{wei2012distributed}, \cite{WeiOzdaglar2013arX}, \cite{shi2014linear}, \cite{chang2015multi}. 
In the augmented space, due to separability of DOP \eqref{miniBis}, monotonicity of the full gradient is automatically maintained. Different from that, in a game setup, individual joint convexity is too restrictive, unless the game is separable to start with, which is a trivial case. Rather, a typical assumption is  individual \emph{partial} convexity of each cost, and  monotonicity  (strict, strong) of the pseudo-gradient. 
Furthermore,  in the augmented space of actions and estimates, due to the inherent coupling  in \eqref{mini_00}, monotonicity is not necessarily maintained. {We note that in ADMM for DOP \eqref{miniBis}, relaxation from strongly convex to merely convex $f_i$ involves  either  an exact augmented minimization at each step, (hence $f_i$ needs to have inexpensive proximal-operator), \cite{WeiOzdaglar2013arX}, or an extra quadratic augmentation, for a sufficiently large step-size,   \cite{HongLuoRaza2016}.  In ADMM for game \eqref{mini_00}, we use a similar idea as \cite{HongLuoRaza2016}  for the extended pseudo-gradient (which may be non-monotone in the augmented space). However, unlike \cite{HongLuoRaza2016},  since cost functions $J_i$ are not separable,  our analysis is novel, based on the decomposition of the augmented space into the consensus subspace and its orthogonal complement.    
 }
Moreover, our proof techniques rely on the pseudo-gradient, since we cannot employ a common Lagrangian as in DOP. Another consequence is that, unlike \cite{WeiOzdaglar2013arX}, we are  able to characterize the rate of convergence only in terms of the residuals relative to NE optimality. \vspace{-0.25cm}
}
The paper is organized as follows. The problem statement and assumptions are provided in Section~2. In Section~3, an inexact ADMM-type algorithm is developed.  Its convergence  to Nash equilibrium is analyzed in Section~4. Simulation results are given in Section~5 and conclusions in Section~6.\vspace{-0.25cm}
\subsection{Notations and Background}

\vspace{-0.255cm}
A vector $x\in\mathbb{R}^n$ may be represented in equivalent ways as $x=[x_1,\ldots,x_n]^T=[x_1;\ldots;x_n]$, or $x=[x_i]_{i =1,\ldots,n}$, or $x=(x_i,x_{-i})$. All  vectors are assumed to be column vectors. Given a vector  $x\in\mathbb{R}^n$,  $x^T$ denotes its transpose and $\|x\| = \sqrt{x^Tx}$ denotes its Euclidean norm. Given a symmetric $n\times n$ matrix $A$, $\|x\|_A$ denotes the weighted norm $\|x\|_A:=\sqrt{x^TAx}$. Denote $\one_n=[1,1,\ldots,1]^T\in\mathbb{R}^n$ and $\zero_n=[0,0,\ldots,0]^T\in\mathbb{R}^n$. Let $e_i$ denote the $N\times 1$, $i$-th  unit vector in $\mathbb{R}^n$, whose $i$-th element is 1 and the rest are 0.  Given a vector $x \in\mathbb{R}^n$, $diag([x_i]_{i= 1,\ldots,n})$ denotes the $n\times n$ diagonal  matrix with $x_1,\ldots,x_n$ on the diagonal. Similarly, given $A_i$, $i=1,\dots,n$ as $(p\times q)$ matrices,   $ diag([A_i]_{i= 1,\ldots,n})$ denotes  the $(np\times nq) $ block-diagonal matrix with $A_i$ on the diagonal. $I_n$ denotes the identity matrix in $\mathbb{R}^{n\times n}$. The Kronecker product of matrices A and B is denoted by $A\otimes B$. Given  matrices $A,B\in\mathbb{R}^{n\times n}$, $A\succ0$ $(A\succeq0)$  denotes that $A$ is positive (semi-)definite. 
$\lambda_{\min}(A)$ and $\lambda_{\max}(A)$ represent the minimum and the maximum eigenvalue of matrix $A$, respectively. {Given a $m \times n$ matrix $A$,  $\|A\|_2$ denotes its induced 2-norm, $\|A\|_2 = \sigma_{\max} (A)$, $\sigma_{\max}(A)= \sqrt{\lambda_{\max}(A^TA)}$}. For any $a, b\in\mathbb{R}^n$ and $\rho >0$, \vspace{-0.4cm}
\begin{equation}\label{inequality}
-\frac{1}{2\rho}\|a\|^2-\frac{\rho}{2}\|b\|^2\leq a^Tb\leq\frac{1}{2\rho}\|a\|^2+\frac{\rho}{2}\|b\|^2.
\end{equation}\vspace{-0.4cm}
For every $a,b,c\in\mathbb{R}^n$ and $n \times n$ matrix $A\succeq 0$,
\begin{equation}\label{simple_math}
\hspace{-0.15cm}(a-b)^TA(a-c)=\frac{1}{2}\|a-c\|_A^2+\frac{1}{2}\|a-b\|_A^2-\frac{1}{2}\|b-c\|_A^2.
\end{equation}
%
%
Given a set $\Omega\in\mathbb{R}^n$, $|\Omega|$ denotes the cardinality of $\Omega$. The Euclidean projection of $x\in\mathbb{R}^n$ onto  $\Omega\subset\mathbb{R}^n$ is denoted by $T_\Omega\{x\}$. Denote by $\text{prox}_{g}^{a}$ the proximal operator for function $g$ with a constant $a$, defined as follows: \vspace{-0.4cm}
\begin{equation}\label{prox_def}
\text{prox}_{g}^{a}\{s\}:=\text{arg }\min_{x}\{g(x)+\frac{a}{2}\|x-s\|^2\}. \vspace{-0.255cm}
\end{equation}
Let $\mathcal{I}_{\Omega}(x):=\begin{cases}0&\text{if }x\in\Omega\\\infty&\text{otherwise}\end{cases}$, be  the indicator function of a closed convex set $\Omega \in\mathbb{R}^n$. Then, $
\text{prox}_{\mathcal{I}_{\Omega}}^{a}\{\cdot\}=T_{\Omega}\{\cdot\}$. 
{Let  $\partial \mathcal{I}_{\Omega}(x)$ denote the subdifferential of $\mathcal{I}_{\Omega}$ at $x$, i.e., the set of all subgradients of $\mathcal{I}_{\Omega}$ at $x$ .  
Then $\partial \mathcal{I}_{\Omega}(x)$  is a convex set and $\partial \mathcal{I}_{\Omega}(x) = N_{\Omega}(x)$, where 
 $  N_{\Omega}(x) = \{  y \in\mathbb{R}^n | \, y^T ( x'-x) \leq 0, \forall x'\in \Omega \}$  is  the normal cone to $\Omega$ at $x$. Moreover, 
 $(y_1-y_2)^T(x_1-x_2) \geq 0$, $\forall x_1, x_2 \in \Omega$, $\forall y_1 \in \partial \mathcal{I}_{\Omega}(x_1)$, $\forall y_2 \in \partial \mathcal{I}_{\Omega}(x_2)$,  \cite{Rockafellar1970}.
 }

\vspace{-0.255cm}
For an undirected graph $G(V,E)$, we denote by:
\begin{itemize}\vspace{-0.25cm}
	\item $V$: Set of vertices in $G$,
	\item $E\subseteq V\times V$: Set of all edges in $G$. $(i,j)\in E$ if and only if $i$ and $j$ are connected by an edge,
	\item $N_i:=\{j\in V|(i,j)\in E\}$: Set of neighbours of $i$ in $G$,
	\item $A:=[a_{ij}]_{i,j\in V}$: Adjacency matrix of $G$ where $a_{ij}=1$ if $(i,j)\in E$ and $a_{ij}=0$ otherwise,
	\item $\mathbb{D}:=\text{diag}(|N_1|,\ldots,|N_N|)$: Degree matrix of $G$,  $d = trace(\mathbb{D}) = \sum_{i=1}^N |N_i|$, $d^*=\max_{i} \{ |N_i | \}$, 
	\item $L:\mathbb{D}-A$: Laplacian matrix of $G$,
	\item $L_N:=\mathbb{D}^{-\frac{1}{2}}L\mathbb{D}^{-\frac{1}{2}}$: Normalized Laplacian of $G$ if $G$ has no isolated vertex.   \vspace{-0.255cm}
\end{itemize}
The following hold for  a graph $G$ with no isolated vertex: $\mathbb{D}\succ 0$, $\lambda_{\max}(L_N)\leq2$ and  $2\mathbb{D}-L = \mathbb{D}^{\frac{1}{2}}(2I-L_{N})\mathbb{D}^{\frac{1}{2}} \succeq0$, \cite{chung1997spectral}. 
%
For a connected, undirected graph $G$ with $n$ vertices, $L\succeq 0$, with 0 a simple eigenvalue and $L\one_n=\zero_n$, $\one_n^TL=\zero_n^T$.  The other eigenvalues of $L$ are positive, with minimum one  $\lambda_2(L)$ and maximum $\lambda_{\max}(L)\leq 2d^*$. { 
Also, $Ker(L) = span \{ \mathbf{1}_n \}$, $Ker(L)^{\perp} = \{ x |  \mathbf{1}^T_n x = 0\}$,  $x^T L x \geq \lambda_2(L) \|x\|^2 >0$,  $\forall x \neq 0$, $x \in Ker(L)^{\perp}$.
 }
\
 \vspace{-0.255cm}
 \section{Problem Statement}\label{problem_statement}

 \vspace{-0.255cm}
Consider a networked game with $N$ players,  defined with the following parameters and denoted by $\mathcal{G}(V,\Omega_i,J_i)$:\vspace{-0.25cm}
\begin{itemize}
	\item $V=\{1,\ldots,N\}$: Set of all players,
	\item $\Omega_i\subset\mathbb{R}$: Action set of player $i$, $\forall i\in V$,
	\item $\Omega=\prod_{i\in V}\Omega_i\subset\mathbb{R}^N$: Action set of all players, where $\prod$ denotes the Cartesian product, 
	\item $J_i:\mathbb{R}^N \rightarrow \mathbb{R}$: Cost function of player $i$, $\forall i\in V$.  
\end{itemize} \vspace{-0.255cm} 
Players' actions are denoted as follows:\vspace{-0.25cm}
\begin{itemize}
	\item $x_i\in\Omega_i$: Player $i$'s action, $\forall i\in V$,
	\item $x_{-i}\in\Omega_{-i}:=\prod_{j\in V\backslash\{i\}}\Omega_j$: All players' actions except player $i$'s,
	\item $x=(x_i,x_{-i})\in\Omega$: All players actions.
\end{itemize}\vspace{-0.25cm}
The game is played such that for a given $x_{-i}\in \Omega_{-i}$, every player $i\in V$ aims to minimize its own cost function { \eqref{mini_00}} with respect to (w.r.t.) $x_i$.  
Note that each player's optimal action is dependent on the other players' decisions. A Nash equilibrium (NE) lies at the intersection of solutions to the set of problems \eqref{mini_00} (fixed-point of best-response map),  such that no player can reduce its cost by unilaterally deviating from its action.  \vspace{-0.3cm} 
\begin{definition}\label{Nash_def}
	Consider an $N$-player game $\mathcal{G}(V,\Omega_i,J_i)$. 
	An action profile (vector) $x^*=(x_i^*,x_{-i}^*)\in\Omega$ is called a Nash equilibrium (NE) of this game if \vspace{-0.25cm}
$$	J_i(x_i^*,{x_{-i}^{*}})\leq J_i(x_{i},{x_{-i}^{*}})\quad\forall x_i\in \Omega_i,\,\,\forall i\in V.$$ \vspace{-0.45cm}
\end{definition} \vspace{-0.45cm}
 We state a few assumptions for the existence of a NE, 
\cite{Nash1950,Debreu1952,Glicksberg1952}.  \vspace{-0.255cm}
\begin{assumption}
	\label{assump}
	For every $i\in V$, \vspace{-0.25cm}
	\begin{itemize}
		\item $\Omega_i\subset\mathbb{R}$ is non-empty, compact and convex,
		\item $J_i(x_i,x_{-i})$ is $C^1$ and convex in $x_i$, for every $x_{-i}$, and jointly continuous in $x$.
	\end{itemize}
\end{assumption} \vspace{-0.255cm}
Let   $\nabla_iJ_i(x) = \frac{\partial J_i}{\partial x_i}(x_i, x_{-i})$  be the partial gradient of $J_i$ w.r.t. $x_i$ 
and   $F :{\mathbb{R}^N}\rightarrow \mathbb{R}^N$ be the pseudo-gradient of game \eqref{mini_00} defined by  \vspace{-0.255cm}  
\begin{equation}\label{pseudoDef}
F(x):=[\nabla_iJ_i(x)]_{i\in V} 
\vspace{-0.255cm}
\end{equation}
Let $\partial\mathcal{I}_{\Omega_i}(x_i)$ be the {subdifferential} of $\mathcal{I}_{\Omega_i}$ at $x_i$  
and $G(x):=[\partial\mathcal{I}_{\Omega_i}(x_i)]_{i\in V}$, where {   $\partial \mathcal{I}_{\Omega_i}(x_i) =N_{\Omega_i}(x_i)$. }
{Let $N_\Omega(x) = \prod_{i\in V} N_{\Omega_i} (x_i)$ be the normal cone to $\Omega$ 
at $x$. 
}

\vspace{-0.255cm}
{    
A Nash equilibrium (NE)  $x^* = (x_i^*,{x_{-i}^{*}}) \in \Omega$ satisfies the variational inequality, (cf. Proposition 1.4.2, \cite{Facchi24})\vspace{-0.25cm}
$$(x-x^*)^T F(x^*) \geq 0, \forall x \in \Omega,\vspace{-0.25cm}$$i.e.,
$-F(x^*) \in N_\Omega(x^*) $, or \vspace{-0.25cm}
$$- \nabla_iJ_i(x^*) \in N_{\Omega_i}(x_i^*),\,\forall i\in V,
$$}Thus, with {   $\partial \mathcal{I}_{\Omega_i}(x_i) =N_{\Omega_i}(x_i)$ }, 
 a NE $x^*= (x_i^*,{x_{-i}^{*}})$ of \eqref{mini_00} can be characterized  by  
\vspace{-0.255cm}
\begin{equation}\label{Nime_Nash}
0 \,  { \in}  \, \nabla_iJ_i(x^*)+\partial\mathcal{I}_{\Omega_i}(x_i^*), \,\,\forall i\in V,
\end{equation}
Then, with \eqref{pseudoDef}, \eqref{Nime_Nash} can be written in compact form as, \vspace{-0.25cm}
\begin{equation}\label{Nash_asle_Kari}
\mathbf{0}_N  \, {\in} \, F(x^*)+G(x^*).  \vspace{-0.255cm}
\end{equation}
Typically another assumption such as monotonicity (strict, strong) of the pseudo-gradient vector $F$, \eqref{pseudoDef}, is used to show that projected-gradient type algorithms converge to $x^*$,  \cite{Facchi24}. 
 \vspace{-0.355cm}
{
\begin{assumption}\label{strgmon_Fassump}
 	The pseudo-gradient $F$ is strongly monotone and Lipschitz continuous: there exists $\mu>0$, $\theta_0 >0$ such that for any $x$ and $y$, 
	$\langle x-y, F(x) - F(y) \rangle \geq \mu  \| x- y\|^2$  and   
	$\| F(x) - F(y) \| \leq \theta_0  \| x- y\|$.	   
\end{assumption} \vspace{-0.25cm}
 Under Assumptions \ref{assump} and  \ref{strgmon_Fassump},  the game has  a unique NE $x^*$ (Theorem 2.3.3 in \cite{Facchi24}). 
 Strong monotonicity  of $F$ is a standard assumption under which convergence of projected-gradient type algorithms is guaranteed with fixed step-sizes, (Theorem 12.1.2 in \cite{Facchi24}).
}

 \vspace{-0.25cm}
We assume that the cost function $J_i$ and the action set $\Omega{_i}$ information are available to player $i$ only, { hence an incomplete-information game, \cite{SB87}. 
The challenge is that each optimization problem in \eqref{mini_00} is dependent on the solution of the other simultaneous problems.  In the classical setting of a game with incomplete  information but  perfect monitoring, each agent can observe the actions of all other players, $x_{-i}$,   \cite{Facchi24}, \cite{Scutari2014}, \cite{frihauf2012nash}, \cite{shanbhag}.   In this paper, we consider a game with \emph{ incomplete information and imperfect monitoring}, where no player is able to observe the actions of all  others. We refer to this as a game with \emph{partial-decision information}. To compensate for the lack of global decision information, players  exchange some information in order to update their actions. We assume players can communicate only locally, with their neighbours. }
 An undirected \emph{communication graph} $G_c(V,E)$ is then defined with no isolated vertex. Let denote the set of neighbours of player $i$ in $G_c$ by $N_i$. Let also denote $\mathbb{D}$ and $L$ be the degree and Laplacian matrices associated to $G_c$, respectively. 
The following assumption is used.  \vspace{-0.255cm}
\begin{assumption}\label{connectivity}
	$G_c$ is  undirected and connected.
\end{assumption}  \vspace{-0.255cm}
We assume that players maintain estimates of the other players' actions and share them with their neighbours in order to update their estimates. {In  an engineering application, this step could be prescribed, e.g. as a network protocol in peer-to-peer networks or ad-hoc networks.}  Our goal is to develop a distributed algorithm for computing the NE of $\mathcal{G}(V,\Omega_i,J_i)$ with fixed step-sizes, while using only partial-decision information over the communication graph $G_c(V,E)$. \vspace{-0.3cm}  

{This is a topic of recent interest given the the proliferation of engineering networked applications requiring distributed protocols that operate under partial, local information (e.g. ad-hoc networks or smart-grid networks), and the deregulation of global markets. We provide next a motivating example,  inspired by the one in \cite{Jayash}, \cite{shanbhag}. }   
\vspace{-0.25cm}
\begin{example}\label{Ex_quadratic}{(Nash-Cournot Game Over a Network). 
Consider a networked Nash-Cournot game,  as in \cite{Jayash}, \cite{Bimpikis2014} between a set of $N$ firms (players/agents)  involved in the production of a homogeneous commodity  that compete over $m$ markets, $M_1,\cdots,M_m$ (Figure \ref{fig_network_cournot_game}). Firm $i$, $i \in V$  participates in  $n_i \leq m$ markets with $x_i$ commodity amount that it supplies  to each of its markets. Hence player $i$'s action (strategy)  is  its commodity amount per market   $ x_i \in \mathbb{R}$. Its total amount it produces is $x_{i,T}=n_i x_i$, assumed to be limited to $\Omega_T$, so that $ x_i \in \Omega_i \subset \mathbb{R}$, $\Omega_i = \Omega_T/n_i$. Firm $i$ has a local vector $A_i\in \mathbb{R}^{m}$  (with elements $1$ or $0$) that specifies which markets it participates in. The $k$-th element of $A_i$ is $1$ if and only if player $i$ delivers $x_i$ amount to market $M_k$.  Therefore,  $A_1,\cdots,A_N$ can be used to specify a  bipartite graph that represents the connections between firms and markets (see Figure \ref{fig_network_cournot_game}). Denote  $ x=[x_i]_{i  \in V}\in \mathbb{R}^N$,  and   $A=[A_1,\cdots,A_N] \in \mathbb{R}^{m\times N}$. Then $A \, x\in \mathbb{R}^{m}= \sum_{i=1}^N A_i x_i$ is the total product supply to all markets, given the action profile $x$ of all firms. Suppose that $P: \mathbb{R}^m \rightarrow \mathbb{R}^m$ is a price vector function that maps the total supply of each market to the corresponding market's price. Each firm $i$ has a production cost  $c_i(x_{i,T}): \Omega_T \rightarrow \mathbb{R}$, function of its  total production amount. Then the local objective function of firm (player) $i$ is $J_i(x_i,x_{-i})= c_i(x_{i,T})-P^T(Ax)A_ix_i$. Overall, given the other firms' profile $x_{-i}$, each firm needs to solve the following optimization problem, }\vspace{-0.3cm}
\begin{equation}\label{network_cournot_game}
\begin{aligned}
& \underset{x_i \in \Omega_i}{\text{minimize}}
& & c_i(x_{i,T})-P^T(Ax)A_ix_i \\
\end{aligned}
\end{equation}

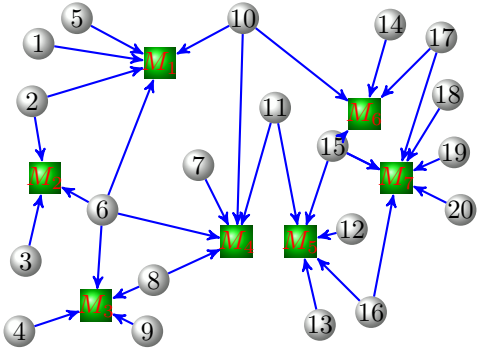
\begin{figure}
\begin{center}
\begin{tikzpicture}[->,>=stealth',shorten >=0.3pt,auto,node distance=2.1cm,thick,
  rect node/.style={rectangle, ball color={rgb:red,0;green,20;yellow,0},font=\sffamily,inner sep=1pt,outer sep=0pt,minimum size=12pt},
  wave/.style={decorate,decoration={snake,post length=0.1mm,amplitude=0.5mm,segment length=3mm},thick},
  main node/.style={shape=circle, ball color=green!1,text=black,inner sep=1pt,outer sep=0pt,minimum size=12pt},scale=0.85]


  \foreach \place/\i in {{(-3.1,3.1)/1},{(-3.2,2.2)/2},
  {(-3.3,-0.3)/3},
  {(-3.4,-1.4)/4},
  {(-2.5,3.5)/5},
  {(-2.1,0.5)/6},
  {(-0.6,1.2)/7},
  {(-1.3,-0.6)/8},
  {(-1.4,-1.4)/9},
  {(0.1,3.5)/10},
  {(0.6,2.1)/11},
  {(1.8,0.2)/12},
  {(1.3,-1.3)/13},
  {(2.4,3.4)/14},
  {(1.5,1.5)/15},
  {(2.1,-1.1)/16},
  {(3.2,3.2)/17},
  {(3.3,2.3)/18},
  {(3.4,1.4)/19},
  {(3.5,0.5)/20}}
    \node[main node] (a\i) at \place {};

      \node at (-3.1,3.1){\rm \color{black}{$1$}};
      \node at (-3.2,2.2){\rm \color{black}{$2$}};
      \node at (-3.3,-0.3){\rm \color{black}{$3$}};
      \node at (-3.4,-1.4){\rm \color{black}{$4$}};
      \node at (-2.5,3.5){\rm \color{black}{$5$}};
      \node at (-2.1,0.5){\rm \color{black}{$6$}};
      \node at (-0.6,1.2){\rm \color{black}{$7$}};
      \node at (-1.3,-0.6){\rm \color{black}{$8$}};
      \node at (-1.4,-1.4){\rm \color{black}{$9$}};
      \node at (0.1,3.5){\rm \color{black}{${10}$}};
      \node at (0.6,2.1){\rm \color{black}{${11}$}};
      \node at (1.8,0.2){\rm \color{black}{${12}$}};
      \node at (1.3,-1.3){\rm \color{black}{${13}$}};
      \node at(2.4,3.4) {\rm \color{black}{${14}$}};
      \node at (1.5,1.5){\rm \color{black}{${15}$}};
      \node at (2.1,-1.1){\rm \color{black}{${16}$}};
       \node at (3.2,3.2){\rm \color{black}{${17}$}};
      \node at (3.3,2.3){\rm \color{black}{${18}$}};
      \node at(3.4,1.4) {\rm \color{black}{${19}$}};
      \node at (3.5,0.5){\rm \color{black}{${20}$}};

  \foreach \place/\x in {{(-1.2,2.8)/1},{(-3,1)/2},{(-2.2,-1)/3},
    {(0,0)/4}, {(1,0)/5}, {(2,2)/6}, {(2.5,1)/7}}
  \node[rect node] (b\x) at \place {};

      \node at (-1.2,2.8){\rm \color{red}{$M_1$}};
      \node at (-3,1){\rm \color{red}{$M_2$}};
      \node at (-2.2,-1){\rm \color{red}{$M_3$}};
      \node at (0,0){\rm \color{red}{$M_4$}};
      \node at (1,0){\rm \color{red}{$M_5$}};
      \node at (2,2){\rm \color{red}{$M_6$}};
      \node at (2.5,1){\rm \color{red}{$M_7$}};

         \path[->,blue,thick]               (a1) edge (b1);

         \path[->,blue,thick]               (a2) edge (b1);
         \path[->,blue,thick]               (a2) edge (b2);

         \path[->,blue,thick]               (a3) edge (b2);

         \path[->,blue,thick]               (a4) edge (b3);

         \path[->,blue,thick]               (a5) edge (b1);

         \path[->,blue,thick]               (a6) edge (b1);
         \path[->,blue,thick]               (a6) edge (b2);
         \path[->,blue,thick]               (a6) edge (b3);
         \path[->,blue,thick]               (a6) edge (b4);

         \path[->,blue,thick]               (a7) edge (b4);

         \path[->,blue,thick]               (a8) edge (b3);
         \path[->,blue,thick]               (a8) edge (b4);

         \path[->,blue,thick]               (a9) edge (b3);

         \path[->,blue,thick]               (a10) edge (b1);
         \path[->,blue,thick]               (a10) edge (b4);
         \path[->,blue,thick]               (a10) edge (b6);

         \path[->,blue,thick]               (a11) edge (b4);
         \path[->,blue,thick]               (a11) edge (b5);

         \path[->,blue,thick]               (a12) edge (b5);

         \path[->,blue,thick]               (a13) edge (b5);

         \path[->,blue,thick]               (a14) edge (b6);

         \path[->,blue,thick]               (a15) edge (b5);
         \path[->,blue,thick]               (a15) edge (b6);
         \path[->,blue,thick]               (a15) edge (b7);

         \path[->,blue,thick]               (a16) edge (b5);
         \path[->,blue,thick]               (a16) edge (b7);

         \path[->,blue,thick]               (a15) edge (b7);

         \path[->,blue,thick]               (a17) edge (b6);
         \path[->,blue,thick]               (a17) edge (b7);

         \path[->,blue,thick]               (a18) edge (b7);
         \path[->,blue,thick]               (a19) edge (b7);
         \path[->,blue,thick]               (a20) edge (b7);
\end{tikzpicture}
\end{center}
\caption{Network Cournot game: An edge from $i$ to $M_k$ on this graph implies that agent/firm $i$ participates in Market $M_k$.}\label{fig_network_cournot_game}
\end{figure}
\begin{figure}
\begin{center}
\begin{tikzpicture}[->,>=stealth',shorten >=0.3pt,auto,node distance=2.1cm,thick,
  rect node/.style={rectangle,ball color=blue!10,font=\sffamily,inner sep=1pt,outer sep=0pt,minimum size=12pt},
  wave/.style={decorate,decoration={snake,post length=0.1mm,amplitude=0.5mm,segment length=3mm},thick},
  main node/.style={shape=circle,ball color=green!1,text=black,inner sep=1pt,outer sep=0pt,minimum size=12pt},scale=0.85]


  \foreach \place/\i in {{(-3,2)/1},{(-2,2)/2},
  {(-1,2)/3},
  {(0,2)/4},
  {(1,2)/5},
  {(2,2)/6},
  {(3,2)/7},
  {(3,1)/8},
  {(3,0)/9},
  {(3,-1)/10},
  {(3,-2)/11},
  {(2,-2)/12},
  {(1,-2)/13},
  {(0,-2)/14},
  {(-1,-2)/15},
  {(-2,-2)/16},
  {(-3,-2)/17},
  {(-3,-1)/18},
  {(-3,-0)/19},
  {(-3,1)/20}}
    \node[main node] (a\i) at \place {};

      \node at (-3,2){\rm \color{black}{$1$}};
      \node at (-2,2){\rm \color{black}{$2$}};
      \node at (-1,2){\rm \color{black}{$3$}};
      \node at (0,2){\rm \color{black}{$4$}};
      \node at (1,2){\rm \color{black}{$5$}};
      \node at (2,2){\rm \color{black}{$6$}};
      \node at (3,2){\rm \color{black}{$7$}};
      \node at (3,1){\rm \color{black}{$8$}};
      \node at (3,0){\rm \color{black}{$9$}};
      \node at (3,-1){\rm \color{black}{${10}$}};
      \node at (3,-2){\rm \color{black}{${11}$}};
      \node at (2,-2){\rm \color{black}{${12}$}};
      \node at (1,-2){\rm \color{black}{${13}$}};
      \node at (0,-2) {\rm \color{black}{${14}$}};
      \node at (-1,-2){\rm \color{black}{${15}$}};
      \node at (-2,-2){\rm \color{black}{${16}$}};
      \node at (-3,-2){\rm \color{black}{${17}$}};
      \node at (-3,-1){\rm \color{black}{${18}$}};
      \node at(-3,0) {\rm \color{black}{${19}$}};
      \node at (-3,1){\rm \color{black}{${20}$}};

        \path[-,blue,thick]               (a1) edge (a2);
            \path[-,blue,thick]               (a2) edge (a3);
         \path[-,blue,thick]               (a3) edge (a4);
         \path[-,blue,thick]               (a4) edge (a5);

         \path[-,blue,thick]               (a5) edge (a6);
         \path[-,blue,thick]               (a6) edge (a7);
         \path[-,blue,thick]               (a7) edge (a8);
         \path[-,blue,thick]               (a8) edge (a9);
           \path[-,blue,thick]               (a9) edge (a10);
         \path[-,blue,thick]               (a10) edge (a11);
             \path[-,blue,thick]             (a11) edge (a12);
            \path[-,blue,thick]               (a12) edge (a13);
            \path[-,blue,thick]               (a13) edge (a14);
            \path[-,blue,thick]               (a14) edge (a15);

            \path[-,blue,thick]               (a15) edge (a16);
            \path[-,blue,thick]               (a16) edge (a17);
            \path[-,blue,thick]               (a17) edge (a18);

            \path[-,blue,thick]               (a18) edge (a19);
            \path[-,blue,thick]               (a19) edge (a20);
             \path[-,blue,thick]               (a1) edge (a20);

              \path[-,blue,thick]               (a2) edge (a15);
             \path[-,blue,thick]               (a6) edge (a13);

\end{tikzpicture}
\end{center}
\caption{Communication graph $G_c$: Firms $i$ and $j$ are able to exchange their local $x^i$ and $x^j$ if there exists an edge between them on this graph. }\label{fig_G_c_graph}
\end{figure}
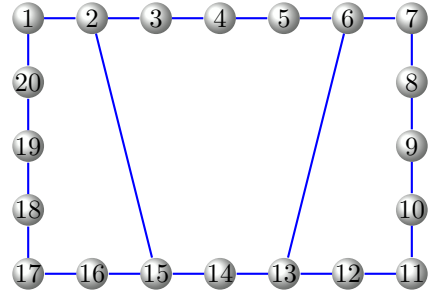
\vspace{-0.5cm}
{In the classical,  centralized-information setting each firm is assumed to have instantaneous access to the others' actions $x_{-i}$. This may be impractical in a large network of geographically distributed firms, \cite{Bimpikis2014}. For example, consider that Figure \ref{fig_network_cournot_game} depicts a group of $N=20$ firms located on different continents that participate in $m=7$ markets, with no centralized communication system between them.  Since players  are unable to directly observe the actions of all others, they  engage in local, non-strategic information exchange, to mitigate their lack of global, centralized information. Firms/players  
may communicate with a local subset of neighbouring firms  via some underlying \textit{communication infrastructure}, hence a distributed partial-information setting.  A communication network formed $G_c$ between the firms prescribes how they communicate locally their production decision,  \cite{Jayash}.  In this situation, the communication network is formed by the players who are viewed as the nodes in the network.  In this example, we consider that most of the communication is between firms on the same continent, with one or at most two firms in each continent having a direct connection to another firm on another continent.  One such instance of the communication network  $G_c$ is  shown in Figure \ref{fig_G_c_graph}. 
Firms $i$ and $j$ are able to exchange their local $x^i$ and $x^j$ if there exists an edge between them on this graph.
Various other topologies can be considered, with different connectivity.}  
\end{example} \vspace{-0.255cm}
{A Nash-Cournot game  \eqref{network_cournot_game} can describe other practical decision problems in engineering networks, 
 for example rate control games in communication networks (\cite{yin2011nash}) 
and demand-response games in smart-grid networks (\cite{YeHu2017}, \cite{peng2}). 
Another  example, of a rate control game over a wireless ad-hoc network, is  presented in Section \ref{sec:simulations}.
}

\section{Distributed Inexact-ADMM Algorithm}\label{asynch}
In order to find the NE of  $\mathcal{G}(V,\Omega_i,J_i)$~\eqref{mini_00} over the communication graph $G_c(V,E)$, we propose a distributed ADMM-type algorithm with fixed step-sizes  (\cite{bertsekas1999parallel}, p. 255), based on   
introducing local estimates of players' actions and using (virtual) constraints for estimate consensus. \vspace{-0.3cm} 

We define a few notations for players' estimates.\vspace{-0.25cm}
\begin{itemize}
	\item $x_{-i}^i\in \mathbb{R}^{N-1}$: Player $i$'s estimate of all other players' actions, 
	\item $x_i^i\in\Omega_i \subset \mathbb{R}$: Player $i$'s estimate of its action which is indeed its own action, i.e., $x_i^i=x_i$ for $i\in V$, 
	\item $x^i=(x_i^i,x_{-i}^i)\in \mathbb{R}^{N}$: Player $i$'s estimate of all players' actions (state variable), 
	\item $\mathbf{x} = [x^i]_{i \in V}$ or $\mathbf{x}=[{x^1};\ldots;{x^N}]\in \mathbb{R}^{N^2}$: Augmented (stacked) vector of all estimates.
\end{itemize}  \vspace{-0.255cm}
Note that all players' actions $x$ can be  represented as $x=[x_i^i]_{i\in V}$.  In steady-state  all local copies $x^i$ should be equal. 
By employing the actions' estimates  $x^i$, $\forall i\in V$, (local copies of $x$),  we can reformulate game \eqref{mini_00} 
as: \vspace{-0.3cm}  
\begin{eqnarray}
\label{mini_1}
&&\hspace{-0.09cm}\begin{cases}
\begin{aligned}
& \underset{x_i^i\in \Omega_i}{\text{minimize}}
& & J_i(x_i^i,x_{-i}^i), \qquad \qquad  \qquad \qquad   \forall i\in V\\  
& \text{subject to}
& &  x^i=x^j,  \, \, \forall j  \text{ s.t. } (i, j), \, (j,i) \in E  \\
\end{aligned}
\end{cases} \vspace{-0.15cm} 
\end{eqnarray}
where    $E$ is the set of edges of $G_c$. For each agent $i$ we consider the constraint that the estimate (local copy)  $x^i$ is equal to  $x^j$  of its neighbours  for all cases where $x^i$ is involved, i.e., player $i$ is either a source or a destination, as specified via ordered pairs $(i, j) \in E$ and  $(j,i) \in E$.   Under Assumption \ref{connectivity} this ensures that 
 $x^1=x^2\ldots =x^N$, hence \eqref{mini_1} recovers \eqref{mini_00}.
The following distributed ADMM-type algorithm is proposed. 

\vspace{-0.15cm}
\begin{algorithm}
	\caption{Inexact-ADMM Algorithm}
	\label{ADMMalgorithm}
	\begin{algorithmic}[1]
		\State \textbf{initialization}  
		
		$x_i^i(0)\in\Omega_i$, $x_{-i}^i(0)\in \mathbb{R} ^{N-1}$, 
		$w^{i}(0)=\textbf{0}_N$ \vspace{0.2cm}
		\For{$k=1,2,\ldots$ }\vspace{0.2cm}
		\For{\hspace{-0.1cm} each player $i\in V$ }\vspace{0.2cm}
		\State \hspace{-0.7cm}  players $i$, $j$ $\forall j\in N_i$		exchange $x^i(k\!-\!1\!)$, $x^j(k\!-\!1\!)$\vspace{0.2cm}
   		\State {\Statex \hspace{-1cm} \hspace{1.08cm} $w^i(k)=w^i(k-1)+c\sum_{j\in N_i}(x^i(k-1)-x^j(k-1))$}
		\State \hspace{-1cm}  $x_i^i(k)=\!T_{\Omega_i}\!\Big\{\!\alpha_i^{-1}\!\Big[ \! -\nabla_iJ_i({x_i^i(k\!-\!1),x^i_{-i}(k\!-\!1)}) \!-\!w_i^i(k)$ 
		\Statex \hspace{0cm}  \hspace{1cm}$+\beta_i x_i^i(k-\!1\!)+\! { \bar{c}}\!\sum_{j\in N_i}\!(x_i^i(k\!-\!1\!)\!+\!x_i^j(k\!-\!1\!))  \Big]\Big\}$\vspace{0.2cm}
		\State \hspace{-1cm} $x_{-i}^i(k)\!=\!\frac{
		\beta_i x_{-i}^i(k-\!1\!)+\!{ \bar{c}}\! \sum_{j\in N_i}\!(x_{-i}^i(k-\!1\!)+x_{-i}^j(k-\!1\!)) - w_{-i}^i(k)}{\alpha_i}$\vspace{0cm}
		\EndFor\vspace{0.2cm}
		\EndFor
	\end{algorithmic}
\end{algorithm}
Here $\alpha_i=\beta_i+2 { \bar{c}} |N_i|$,  { $\bar{c}= c +c_0$},   
 where  $c>0$, $c_0>0$, $\beta_i>0$ are penalty parameters for the augmented Lagrangian and proximal first-order approximation.   
  Besides its estimates, each player $i$ needs to maintain only its dual variables  $w^i$, and needs only the information on its $\Omega_i$ set.   {The derivation of Algorithm \ref{ADMMalgorithm} following an ADMM approach adapted to a game setup is provided in Appendix A.}  
\vspace{-0.25cm}
  
Algorithm \ref{ADMMalgorithm} has a gradient-play structure for a game with estimate consensus constraints: Step 6 uses a projected-gradient  (due to $J_i$) of an  augmented Lagrangian with quadratic-penalty and proximal linearization (inexact ADMM) (see \eqref{x_i^i_ADMM_Pre1} in Appendix A). 
 So in this sense,  in the context  of individual rationality, it can be interpreted as a better-response (as opposed to best-response) strategy, but in the modified game, due to partial-information setting. Step 5 updates the dual variables $w^i$ via a dual ascent. Compared to a gossip-based algorithm,  \cite{Jayash}, \cite{salehisadaghiani2016distributed},  the additional 
 ADMM-induced modifications are  instrumental to prove convergence  under fixed step-sizes. \vspace{-0.255cm}
 
\begin{remark}\label{vs_opt} { 
We compare Algorithm  \ref{ADMMalgorithm} to typical ADMM  DOP algorithms. 
Unlike   \cite{wei2012distributed}, \cite{WeiOzdaglar2013arX}, where a minimization subproblem is solved exactly at each iteration, Algorithm  \ref{ADMMalgorithm} uses a gradient update based on inexact  (proximal) approximation as in \cite{HongLuoRaza2016}.  In contrast to either of these  algorithms, the update of each player decision  $x_i^i$ (Step 6) is indirectly coupled nonlinearly to  the previous decisions of other players, via its estimate $x^i_{-i}$,  due to the intrinsic coupling in  $J_i$.  This is unlike ADMM for DOP where the update of each agent's decision variable $x_i$ is only linearly coupled to the others agents' $x_j$. Moreover, in Algorithm 1 each player  has an extra step for updating its estimate $x^i_{-i}$ (Step 7), in addition to its decision variable $x_i^i$. Because of this coupling the convergence proof (based on the pseudo-gradient) is more involved.} 
\end{remark}\vspace{-0.25cm}


Next we analyze the fixed points  of Algorithm~\ref{ADMMalgorithm}.

\begin{lemma}\label{correctness_lemma_1} 
		Consider that Assumptions~\ref{assump}, \ref{connectivity} hold. Let $\overline{\mathbf{x}}:=[{\overline{x}^i}]_{i\in V}$, 
		$\overline{\mathbf{w}}:=[{\overline{w}^i}]_{i\in V}$ be any fixed points of Algorithm~\ref{ADMMalgorithm}. 
		Then, $\overline{\mathbf{x}} = \mathbf{1}_N \otimes \overline{x}$, where $\overline{x}= x^*$  is  NE of  game \eqref{mini_00}, 
and 
		$(\mathbf{1}^T_N \otimes I_N)\overline{\mathbf{w}} =  \mathbf{0}_{N}$. 
\end{lemma}
\par{\emph{Proof}}. From Step 5 of Algorithm~\ref{ADMMalgorithm},  for all $i\in V$,  \vspace{-0.3cm} 
\begin{equation*}\label{fixed_xi_va_xj}{
\hspace{0cm}\begin{array}{crl}
\overline{w}^i=\overline{w}^i+c\sum_{j\in N_i}(\overline{x}^i-\overline{x}^j)
\end{array} \, \Rightarrow  (L \otimes I_N) \overline{\mathbf{x}}= \mathbf{0}_{N^2} }
\end{equation*}
and by Assumption~\ref{connectivity}, this yields $\overline{\mathbf{x}} = \mathbf{1}_N \otimes \overline{x}$, for some  $\overline{x} =[\overline{x}_i]_{i \in V}\in \mathbb{R}^{N}$, hence  \vspace{-0.3cm} 
\begin{equation}\label{xtildemosavi}\overline{x}^1=\ldots=\overline{x}^N=\overline{x},\end{equation}
Using \eqref{xtildemosavi} in Step 6 of Algorithm~\ref{ADMMalgorithm} yields,
$  (\beta_i+2  { \bar{c}}|N_i|)\overline{x}_i - \alpha_i\overline{x}_i  \in \nabla_iJ_i(\overline{x}) +\partial\mathcal{I}_{\Omega_i}(\overline{x}_i) + \overline{w}_i^i  
$.  
With $\alpha_i = \beta_i + 2 { \bar{c}}|N_i|$ this leads to \vspace{-0.3cm} 
\begin{equation*}\label{fixed_xi_nahayi}
0 \, {\in} \, \nabla_iJ_i(\overline{x})+\partial\mathcal{I}_{\Omega_i}(\overline{x}_i)+\overline{w}_i^i, \, \, \, \forall i\in V \vspace{-0.25cm}
\end{equation*}
Moreover, using \eqref{xtildemosavi} in Step 7 of Algorithm~\ref{ADMMalgorithm}, we obtain $\alpha_i\overline{x}_{-i}^i=\alpha_i\overline{x}_{-i}^i-\overline{w}_{-i}^i
$,  hence  \vspace{-0.25cm}
\begin{equation*}\label{step7nahayi}
\overline{w}^i_{-i}=\zero_{N-1} \, \, \forall i\in V \vspace{-0.25cm}
\end{equation*} 
Combining the previous two relations into a single vector,  yields  \vspace{-0.25cm}
\begin{equation}\label{eqn:barx_NE}
\mathbf{0}_N { \in} \, \nabla_iJ_i(\overline{x})e_i+\partial\mathcal{I}_{\Omega_i}(\overline{x}_i)e_i+\overline{w}^i,  
\, \, \, \forall i\in V \vspace{-0.3cm} 
\end{equation} 
where $e_i $ is the $i$-th unit vector in $ \mathbb{R}^N$. 
{
Summing after $i \in V$, with $F(\overline{x})=[\nabla_iJ_i(\overline{x})]_{i\in V}$, 
yields  \vspace{-0.25cm}
\begin{equation}\label{eqn:barx_w_NE_sum}
\mathbf{0}_N { \in} \, F(\overline{x}) + G(\overline{x})+\sum_{i \in V} \overline{w}^i.  \vspace{-0.255cm}
\end{equation} 
For $w^i(k)$, from Step~5 of Algorithm~\ref{ADMMalgorithm},   it follows that   \vspace{-0.3cm} 
\begin{equation}\label{eqn:wk_sum}
\sum_{i\in V} w^i(k)=\sum_{i\in V} w^i(k-1) + c\sum_{i \in V} \sum_{j\in N_i}(x^i(k-1)-x^j(k-1))  
\end{equation}
Under Assumption  \ref{connectivity} the second term in \eqref{eqn:wk_sum} is zero, hence $\sum_{i\in V} w^i(k)=\sum_{i\in V} w^i(k-1) $ for any $k\geq 1$. With initial conditions $w^i(0)=\textbf{0}_N$ $\forall i\in V$ this implies that  $\sum_{i\in V} w^i(k)= \zero_{N}$ $\forall k\geq 1$, hence $\sum_{i \in V} \overline{w}^i =  \zero_{N}$. 
Thus \eqref{eqn:barx_w_NE_sum} reduces to 
$\mathbf{0}_N{ \in} \, F(\overline{x})+G(\overline{x})$,  
hence by \eqref{Nash_asle_Kari} $\overline{x}=x^*$ is NE  and $\overline{\mathbf{x}} = \mathbf{1}_N \otimes x^*$. }  
Moreover,  \vspace{-0.3cm} 
\begin{equation}\label{barw_0}
(\one_N^T\otimes I_N) \overline{\mathbf{w}} = \sum_{i \in V} \overline{w}^i =  \zero_{N}  \vspace{-0.4cm}
\end{equation} 
$\hfill\blacksquare$

\vspace{-0.2cm}
Recall  the pseudo-gradient \eqref{pseudoDef},  $F(x)=[\nabla_iJ_i(x)]_{i\in V}$. \vspace{-0.3cm} 
\begin{definition}\label{def_bF}
 Let $\mathbf{F}:\mathbb{R}^{N^2}\rightarrow\mathbb{R}^N$,  be the extension of $F$ to  the augmented space,  defined as \vspace{-0.3cm} 
\begin{equation}\label{ExtpseudoDef}
\mathbf{F}(\mathbf{x}):=[\nabla_iJ_i(x^i)]_{i\in V}, \,\, \nabla_i J_i(x^i) = \frac{\partial J_i}{\partial x^i_i}(x^i_i, x^i_{-i}) \vspace{-0.3cm} 
\end{equation}
and called the \emph{extended pseudo-gradient}. Let also  $\mathbf{G}(\mathbf{x}):=[\partial I_{\Omega_i}(x^i_i)]_{i\in V}$, so $\mathbf{G}(\mathbf{x}) = G(x)$.
 \end{definition} \vspace{-0.3cm} 
  Note that $\mathbf{F}$ is the pseudo-gradient evaluated at estimates instead of actual actions,  and 
 $
  \mathbf{F}( \mathbf{1}_N \otimes x)  = F(x)$. 
We write Algorithm~\ref{ADMMalgorithm} in a compact vector form.  \vspace{-0.3cm} 
\begin{proposition}\label{propAlg1_vector}
Let $\mathbf{x} = [x^i]_{i \in V}$, $\mathbf{w} = [w^i]_{i \in V}$.  Then  the updates in Algorithm~\ref{ADMMalgorithm} have the following vector form: \vspace{-0.25cm}
\begin{eqnarray}
&& \mathbf{w}(k) = \mathbf{w}(k-1)+c  \mathbf{L} \mathbf{x}(k-1)  \\  \label{wstep_Algvector}
&& \nonumber \\
&& \mathbf{0}_{N^2} \in  \mathbf{R} \,( \mathbf{\mathbf{F}}( \mathbf{x}(k-1))+\mathbf{\mathbf{G}}( \mathbf{x}(k))) + \mathbf{w}(k) \label{xstep_step67_altVec0}\\
&&\hspace{0.8cm}+\! (( \mathbb{B} \!+ \! 2 {\bar{c}} \mathbb{D}) \!\otimes \! I_N)( \mathbf{x}(k) \!-\!  \mathbf{x}(k \!-\!1))  
\!+\! { \bar{c}} \, \mathbf{L} \mathbf{x}(k\!-\!1),  
\nonumber 
\end{eqnarray}
where $ \mathbf{R} = \text{diag} ([e_i]_{i\in V}) $ is $(N^2 \times N)$ block-diagonal  matrix with unit vectors $e_i$, $\mathbf{L} = L \otimes I_N$,  $\mathbf{F}$,  $\mathbf{G}$ as in \eqref{ExtpseudoDef}, $\mathbb{D}$ is the degree matrix of $G_c$, $L$ is the Laplacian matrix of $G_c$, $\mathbb{B}:=\text{diag}((\beta_i)_{i\in V})$ with $\beta_i>0$, {$\bar{c}=c+c_0$}, $c, c_0>0$. 
\end{proposition} \vspace{-0.25cm}
\par{\emph{Proof}}. 
See Appendix B. \vspace{-0.25cm}

{Properties of $\mathbf{\mathbf{F}}$ and matrices $\mathbf{R}$ and $\mathbf{L}$ play a key role in the following. Matrix $\mathbf{R}^T$ allows action selection from the stacked $\mathbf{x}$: $\mathbf{R}^T  \mathbf{x} =[ e_i^T \, x^i ]_{i \in V} \! =x$, based on  $e_i^T  x^i \! = \! x^i_i \!= \!x_i$. 
Also $ (\mathbf{1}^T_N \otimes I_N)\mathbf{R}  \! = \! I_N$,  $\mathbf{R}^T \mathbf{R} \! =\! I_N$. 
For matrix $\mathbf{L}$,  $Ker (\mathbf{L})\!\! = \!\!Ker (L\! \otimes  \!I_N) \!\!= \!\!span \{ \mathbf{1}_N \! \otimes \! I_N \}$ (the consensus subspace), and $Ker (\mathbf{L})^\perp \! \!= \! \! span (\mathbf{L}) \!= \!Ker \{ \mathbf{1}^T_N\! \otimes \! I_N \}$. }
\vspace{-0.25cm}

{Note that by Lemma \ref{correctness_lemma_1},  $\overline{\mathbf{x}} \!= \!\mathbf{1}_N  \otimes \overline{x} \in Ker (\mathbf{L})$ and   $\overline{\mathbf{w}} \in span (\mathbf{L}) $ (cf. \eqref{barw_0}). Also,  we can write \eqref{eqn:barx_NE}  as,  \vspace{-0.3cm}  
\begin{equation*}
\mathbf{0}_{N^2} { \in} \, [\nabla_iJ_i(\overline{x})e_i]_{i\in V}+[\partial\mathcal{I}_{\Omega_i}(\overline{x}_i)e_i]_{i\in V} + \overline{\mathbf{w}} 
\label{combo_no_stack} \vspace{-0.3cm} 
\end{equation*}
or, with $\mathbf{F}$, $\mathbf{G}$ as in \eqref{ExtpseudoDef}, compactly as,  \vspace{-0.3cm} 
\begin{equation}
\mathbf{0}_{N^2} { \in} \, \mathbf{R} \, ( \mathbf{F} (\overline{\mathbf{x}})+ \mathbf{G} (\overline{\mathbf{x}})) + \overline{\mathbf{w}}. 
\label{combo_no_stackBFG}
\end{equation}
}
\vspace{-1.1cm} 
\section{Convergence Analysis}\label{convergence_proof_Section}
In this section we show global convergence of Algorithm~\ref{ADMMalgorithm}  
under 
the following assumption. \vspace{-0.3cm} 
{
 \begin{assumption}\label{Lipchitz_boldFassump}
 	The extended pseudo-gradient $\mathbf{F}$,  \eqref{ExtpseudoDef},  is Lipschitz continuous: there exists  $\theta>0$ such that for any $\mathbf{x}$ and $\mathbf{y}$,    
$\| \mathbf{F}(\mathbf{x})-\mathbf{F}(\mathbf{y}))\| \leq \theta \| \mathbf{x}- \mathbf{y}\|$.  
	\end{assumption}
} 

 \vspace{-1.cm}
 {
\begin{remark}\label{rem_theta_quadratic} In  the general  case when $\mathbf{F}$ is nonlinear, a sufficient condition for Assumption  \ref{Lipchitz_boldFassump} is that   $\mathbf{F}$ is $C^1$ with bounded Jacobian $\|D \mathbf{F}(\mathbf{x})\|_2$. For quadratic games, Assumption  \ref{Lipchitz_boldFassump} is automatically satisfied ($\mathbf{F}$ is linear). Moreover, as shown below, $\theta =\theta_0$, where $\theta_0$ is the Lipschitz constant of $F$ (Assumption \ref{strgmon_Fassump}). 
Consider the  Nash-Cournot game  in Example \ref{Ex_quadratic}, where firm $i$'s  production cost  is a strongly convex quadratic function in its total production amount,   $c_i(x_{i,T})=  q_i x^2_{i,T}+  b_i x_{i,T}$,  i.e.,  $c_i(x_i)= n^2_i q_i x^2_i+ n_i b_i x_i$, where $q_i >0$,  $b_i\in \mathbb{R}$. Consider that market $M_k$'s price is a linear  function of the total supplied commodity amount, $p_k(x)= \bar{P}_k -  z_k [Ax]_{k}$ (known as a linear inverse demand function) with $\bar{P}_k, z_k>0$. Denote $P=[p_k]_{k=1,m}: \mathbb{R}^N \rightarrow \mathbb{R}^m$,  $\bar{P}=[\bar{P}]_{k=1,m}\in \mathbb{R}^m$, 
$Z=diag ([z_k]_{k=1,m}) \in \mathbb{R}^{m\times m}$, so that $P=\bar{P}-Z Ax$ is the vector price function,  
hence the objective function of player $i$  is   \vspace{-0.3cm}
\begin{eqnarray}
&J_i(x_i,x_{-i}) =c_i(x_i)- (\bar{P}-Z A x)^T A_ix_i  \,\, \text{and } \label{network_cournot_game_function} \\
&\nabla_i J_i(x)= \nabla c_i(x_i) + A_i^T Z A_i x_i-  A_i^T \, (\bar{P}-Z A x), \label{Fi_ex1_0} 
\end{eqnarray}
where $ \nabla c_i(x_i) = 2n^2_i q_i x_i + n_i b_i$. Then, with $x=[x_i]_{i \in V}$, $A=[A_1, \cdots A_N]$, $F(x) = [\nabla_iJ_i(x)]_{i\in V}$ is given by  \vspace{-0.3cm}
$$
F(x)= \nabla c(x)  +  diag ([A_i^T Z A_i]_{i\in V}) \, x+ A^T Z A \, x - A^T\bar{P} 
$$where $\nabla c(x)  = 2 diag ([n^2_i q_i]_{i \in V}) \, x + [n_i b_i]_{i \in V}. $ Combining terms, $F(x)$ can be written 
compactly as,  \vspace{-0.3cm}
\begin{equation}\label{F_ex1}
F(x)= Q \, x + r, \quad \text{ where } \, Q := \Sigma +  A^T Z A, 
\end{equation}
$\Sigma := diag \big ([2n^2_i q_i \!+\!A_i^T Z A_i]_{i \in V} \big)$ and  
$r := [n_i b_i]_{i \in V}  \!-\! A^T\bar{P}$. 
Since   $q_i>0$ and $z_k >0$, $\Sigma \succ 0$ and $ A^T Z A \succeq 0$, hence $Q \succ 0$. 
Thus $F(x)$ is Lipschitz continuous with   $\theta_0 = \|Q\|_2$ where $\|Q\|_2 = \sigma_{\max} (Q)$. 
\\
Using \ \eqref{Fi_ex1_0} with $x$ replaced by $x^i = (x_i,x^i_{-i})$, yields   \vspace{-0.3cm}
$$
\nabla_i J_i(x^i) = \nabla c_i(x_i) + A_i^T Z A_i x_i-  A_i^T \, (\bar{P}-Z A x^i).
$$
Thus, $\mathbf{F}(\mathbf{x}):=[\nabla_iJ_i(x^i)]_{i\in V}$, \eqref{ExtpseudoDef}, is given by   \vspace{-0.3cm}
\begin{eqnarray*}\label{boldF_ex1}
&&\mathbf{F}(\mathbf{x})\! = \!\!\nabla c (x)  \!+ \! diag \!\big ([A_i^T Z A_i]_{i\in V} \!\big) x  \!- \! A^T\bar{P}  \\
&&\hspace{.8cm} \,\! +\!diag\! \big( [A_i^T Z A]_{i \in V} \!\big) \mathbf{x} 
\end{eqnarray*} 
Using $\nabla c (x) $, $\Sigma$ and $r$ defined above, and $\mathbf{R}^T ( I_N \otimes  A^T Z A)=diag \big ( [e_i^T A^T Z A]_{i \in V} \! \big)$ for the last term,  
yields  \vspace{-0.3cm}
$$\mathbf{F}(\mathbf{x}) =  \Sigma  \mathbf{R}^T \mathbf{x} + r+  \mathbf{R}^T (I_N \otimes  A^T Z A) \, \mathbf{x}.
$$
where $x=\mathbb{R}^T \mathbf{x}$ was used. Since $\Sigma$ is block-diagonal,  $\Sigma  \mathbf{R}^T= \mathbf{R}^T \! (I_N \otimes \Sigma)$,  and with $Q$ \eqref{F_ex1} we can write  \vspace{-0.3cm}
 \begin{equation}\label{boldF_ex1_final}
 \mathbf{F}(\mathbf{x})=  \mathbf{\overline{Q}}\,\mathbf{x} + r, \quad \text{ where } \, \mathbf{\overline{Q}}:= \mathbf{R}^T  (I_N \otimes Q).
\end{equation}  
Hence, $\|\mathbf{F}(\mathbf{x}) -\mathbf{F}(\mathbf{y}) \| \leq \|\mathbf{\overline{Q}}\|_2 \| \mathbf{x} - \mathbf{y} \|$, where $\|\mathbf{\overline{Q}}\|_2=\| \mathbf{R}^T (I_N \otimes Q)\|_2 \leq \| \mathbf{R}^T\|_2 \| I_N \otimes Q\|_2 = \|Q\|_2$ since $ \| \mathbf{R}\|_2 =1$ 
and $\| I_N \otimes Q\|_2 = \| Q \|_2$, and thus $\theta =\theta_0$. \end{remark}
}

\begin{remark}
 Note that Assumption~\ref{Lipchitz_boldFassump} on $\mathbf{F}$ is weaker than  strong monotonicity in \cite{YeHu2017}), {or  cocoercivty in \cite{salehiPavelIFAC2017}, \cite{WeiShiACC}}.  In  ADMM  for DOP, cocoercivity of the full gradient is used, which  is automatically satisfied  by   \emph{joint convexity} of each decoupled $f_i$, \eqref{miniBis}, and Lipschitz continuity of its gradient (see \cite{chang2015multi}). In contrast,  in a game, because of coupling  to the others' actions in $J_i$ \eqref{mini_00}  and because of partial convexity (see Assumption~\ref{assump}),  monotonicity of the pseudo-gradient is \emph{not automatically satisfied} when extended to the augmented space. 
{The next result shows how under Assumption~\ref{Lipchitz_boldFassump}, a monotonicity property  can be achieved in the augmented space. 
}
 \end{remark}

 \vspace{-0.5cm} 
{\begin{lemma}\label{boldRF_plus_cL_strgmon_lemma} 
		Consider that Assumptions~\ref{assump}, \ref{strgmon_Fassump}, \ref{connectivity} and \ref{Lipchitz_boldFassump} hold and let  \vspace{-0.25cm}
		\begin{equation}\label{eqboldRF_plus_cL_strgmon_5}
\Psi =  \left [
\begin{array}{ccc}
\frac{\mu}{N} && - \frac{\theta + \theta_0}{2\sqrt{N}}  \\ 		
- \frac{\theta + \theta_0}{2\sqrt{N}}  && c_0  \lambda_2 (L)  - \theta 
\end{array}
\right ]
\end{equation}
Then,  for any $c_0 > c_{\min}$ where  $  c_{\min}  \lambda_2(L) =   \frac{(\theta+\theta_0)^2}{4\mu} +\theta$,  
$\Psi \succ 0$, 	and for  any $ \mathbf{x}$ and any $ \mathbf{y} \in   Ker(\mathbf{L}) $, 
\begin{eqnarray}\label{eqboldRF_plus_cL_strgmon}
&&\hspace{-0.15cm}(\!\mathbf{x}\!-\! \mathbf{y})^T \! \big (\mathbf{R}\mathbf{F}(\mathbf{x})\!-\!\mathbf{R}\mathbf{F}(\mathbf{y})\!+\!c_0 \mathbf{L}(\mathbf{x} \!-\! \mathbf{y}) \big )\!\geq \!\bar{\mu} \|\mathbf{x}\!-\!\mathbf{y}\|^2,
\end{eqnarray}
where $\bar{\mu}:=\lambda_{\min}(\Psi)>0$. 
\end{lemma}}\vspace{-0.25cm}
\par{{\emph{Proof}. See Appendix B. }}\vspace{-0.25cm}
{\begin{remark}\label{rem_c_0_quadratic}
For quadratic games  (see Remark \ref{rem_theta_quadratic}), $\theta=\theta_0=\|Q\|_2$, $F(x)$ is strongly monotone with $\mu = \lambda_{\min}(Q)>0$ (by Assumption  \ref{strgmon_Fassump}), hence  the $c_0$ bound in Lemma \ref{boldRF_plus_cL_strgmon_lemma}  simplifies to $c_0 \lambda_2(L) >  \frac{\|Q\|^2_2}{\mu} + \|Q\|_2$. Furthermore, if $Q$ is symmetric, $\|Q\|_2=\sigma_{\max} (Q) = \lambda_{\max}(Q)$, and this reduces to 
 $c_0 >  \frac{\lambda_{\max}(Q)}{ \lambda_2(L)}( \frac{  \lambda_{\max}(Q)}{ \lambda_{\min}(Q)} + 1) $, showing the trade-off between game properties and communication graph connectivity.  
 \end{remark}}
\vspace{-0.25cm}
{Lemma  \ref{boldRF_plus_cL_strgmon_lemma} shows that,  for a sufficiently large $c_0$, $ \mathbf{R}\mathbf{F}(\mathbf{x}\!)\!+ c_0 \mathbf{L} \mathbf{x}$ is  strongly monotone in a \textit{restricted} set of directions $\mathbf{x} - \mathbf{y}$, where $\mathbf{y} \in Ker(\mathbf{L})$.   We can call $\mathbf{R}\mathbf{F}$ a \textit{pre}-$\mathbf{L}$-monotone function, related to the concept of  \textit{pre-monotone gradient} of a \textit{prox-regular} function, \cite{Rockafellar1996}. 
In fact, being in a restricted set of directions, this is a weaker monotonicity property,  
similar to the notion of \textit{restricted} strong convexity used in DOP and high-dimensional statistical estimation \cite{Zhang2017}, \cite{NegahbanNIPS2010}, \cite{NegahbanNIPS2009}.  
To show this key property, an instrumental step is the decomposition of the augmented  space $\mathbb{R}^{N^2}$ into  the consensus subspace, $Ker(\mathbf{L})$ (where $\mathbf{F}$ is strongly monotone), and its orthogonal complement, $Ker(\mathbf{L})^\perp$ (where $ \mathbf{L}$ is strongly monotone). Based on the Lipschitz continuity of $\mathbf{F}$, for a sufficiently large $c_0$, excess strong monotonicity of $ c_0\mathbf{L}$ can balance the  cross-terms (shortage of monotonicity of $\mathbf{R}\mathbf{F}$) when $\mathbf{x}$ is off the consensus subspace  $Ker(\mathbf{L})$, but not  on  $Ker(\mathbf{L})^\perp$, (see proof in Appendix B). We show next that this property is sufficient to prove convergence of Algorithm~\ref{ADMMalgorithm}.
}

 \begin{theorem}\label{theorem_convergence_Algo1}

	Consider that Assumptions~\ref{assump}, {\ref{strgmon_Fassump},} \ref{connectivity} and { \ref{Lipchitz_boldFassump}} hold.
{Take an arbitrary $c>0$, any  $c_0 > c_{\min}$, with $c_{\min}$ as in Lemma \ref{boldRF_plus_cL_strgmon_lemma},  
and let $\bar{c} = c + c_0$. } If $\beta_i>0$ are chosen such that 
	\begin{equation}\label{condition}
	\lambda_{\min}( {\mathbb{B} + 2\bar{c}}\, \mathbb{D}-c\,L) > {\frac{\bar{\theta}^2}{2 \bar{\mu}} }\vspace{-0.15cm}
	\end{equation}
where $\mathbb{B}:=\text{diag}((\beta_i)_{i\in V})$, 
$\mathbb{D}$ and $L$ are the degree and Laplacian matrices of $G_c$,  {$\bar{\theta}= \theta + 2 c_0 d^*$, $d^*$ is the maximal degree of $G_c$ and $ \bar{\mu}$ is defined in Lemma \ref{boldRF_plus_cL_strgmon_lemma}},  then 	 the sequence $\{x^i(k)\}$ $\forall i\in V$, or   $\{\mathbf{x}(k)\}$,  generated by Algorithm~\ref{ADMMalgorithm} 
converges to 
  NE of game \eqref{mini_00} $x^*$,  (or $\mathbf{1}_N \otimes x^*$). 
\end{theorem}
\par{\emph{Proof}}.
%
{Recall that $\mathbf{G}(\mathbf{x(k)}) \!=\! G(x(k))\! = \![\partial I_{\Omega_i}\!(x_i(k))]_{i\in V}$}. From  
the $\mathbf{x}$-update \eqref{xstep_step67_altVec0} in the vector form of Algorithm~\ref{ADMMalgorithm} (see Proposition \ref{propAlg1_vector}), 
it follows that  {there exists $y \in   G(x(k)) \subset \mathbb{R}^N$ such that, }
  \vspace{-0.3cm} 
\begin{eqnarray*}\label{step67_altVec0}
&&\mathbf{0}_{N^2}=  \mathbf{R} \,( \mathbf{\mathbf{F}}( \mathbf{x}(k\!-\!1))\!+\!{y}) \!+\! \mathbf{w}(k)\nonumber\\ 
&& \hspace{0.8cm} +  ((\mathbb{B} \!+\!2 { \bar{c}}\,  \mathbb{D} )\! \otimes \! I_N)( \mathbf{x}(k) \!- \! \mathbf{x}(k\!-\!1))  
+  { \bar{c}}\,   \mathbf{L} \mathbf{x}(k \!-\!1)  
 \nonumber
\end{eqnarray*}
or, substituting    $\mathbf{w}(k)$ by $ \mathbf{w}(k \!+\!1) \!-\! c\,   \mathbf{L} \mathbf{x}(k) $ (cf. \eqref{wstep_Algvector}),  
 \vspace{-0.3cm} 
\begin{eqnarray*}
\mathbf{0}_{N^2} &&=  \mathbf{R} \, \mathbf{\mathbf{F}}( \mathbf{x}(k-1))+ \mathbf{R} \, y  + \mathbf{w}(k+1)- c \,  \mathbf{L} \mathbf{x}(k) \nonumber\\ 
&& 
+  ((\mathbb{B}\! +\!2 { \bar{c}}\,  \mathbb{D} )\! \otimes I_N)( \mathbf{x}(k)\!- \! \mathbf{x}(k\!-\!1)) \! +\! { \bar{c}} \,  \mathbf{L} \mathbf{x}(k\!-\!1).
 \nonumber
\end{eqnarray*}
Using { $\bar{c} = c +c_0$} and  $H:=\mathbb{B} +2 { \bar{c}}\, \mathbb{D}- c \,L$ this yields   \vspace{-0.3cm}
\begin{eqnarray}\label{step67_altVec}
 \mathbf{0}_{N^2} &&= \mathbf{R} \,\mathbf{\mathbf{F}}( \mathbf{x}(k-1))+ \mathbf{R} \, y + \mathbf{w}(k+1)\\ 
&&
 + (H \otimes I_N) \!( \mathbf{x}(k) \!-\!  \mathbf{x}(k\!-\!1\!)) \! +\!{ c_0}\,   \mathbf{L} \mathbf{x}(k\!-\!1\!).\nonumber
\end{eqnarray}
Note that  $H:=\mathbb{B} + c \,(2 \mathbb{D}- L)+2c_0 \mathbb{D}$ hence $H\otimes I_N {\succ} 0$, since $\mathbb{B} {\succ} 0$ and $c\, (2\mathbb{D}-L)\succeq0$.   \vspace{-0.3cm} 


{
Consider  NE $x^*$ and let $\mathbf{x}^* = \mathbf{1}_N \otimes x^*$. Then,  from \eqref{Nash_asle_Kari}, with  $ \mathbf{F}( \mathbf{1}_N \otimes x^*)= F(x^*)$, 
and the definition of  $\mathbf{R}$, it follows that $ \mathbf{0}_{N^2} { \in} \,  \mathbf{R} \mathbf{F}(\mathbf{x}^*) +  \mathbf{R} \mathbf{G}(\mathbf{x}^*)$, i.e., 
there exists  $y^* \in \mathbf{G}(\mathbf{x^*})$ such that  
$ \mathbf{R} \, \mathbf{F} (\mathbf{x}^*)+ \mathbf{R} \, {y^*}  =\zero_{N^2}. 
$
 Moreover, there exists $ \mathbf{q}^* \in  \mathbb{R}^{N^2}$ (any  $ \mathbf{q}^* \in Ker(\mathbf{L})$ with structure $ \mathbf{1}_N \otimes \widetilde{q} $, for $\widetilde{q} \in \mathbb{R}^N$) such that  \vspace{-0.3cm} 
\begin{equation}
\zero_{N^{2}} =\mathbf{R} \, \mathbf{F} (\mathbf{x}^*)+ \mathbf{R} \, {y^*}  +c \mathbf{L} \mathbf{q}^*. 
\label{combo_no_stackBFGVec_ZERO_q}  \vspace{-0.3cm} 
\end{equation}
}
Subtracting \eqref{combo_no_stackBFGVec_ZERO_q} from \eqref{step67_altVec}, multiplying  by $(\mathbf{x}(k) - \mathbf{x}^*)^T$ and using $\mathbf{R}^T \, \mathbf{x} = x$, yields  \vspace{-0.3cm}
\begin{eqnarray}\label{optim_equation_Nash_combo_multiplic_Revised_SumVec}
&&\hspace{-0.15cm}\big(\mathbf{R} \mathbf{F}(\mathbf{x}(k-1))- \mathbf{R} \mathbf{F}(\mathbf{x}^*)\big)^T(\mathbf{x}(k)-\mathbf{x}^*)\nonumber\\
&&\hspace{-0.15cm}+\big({y-y^*}\big)^T(x(k)-x^*)\nonumber\\  
&&\hspace{-0.15cm}+(\mathbf{w}(k+1)- {c  \mathbf{L} \mathbf{q}^*})^T(\mathbf{x}(k)-\mathbf{x}^*)\\ 
&&\hspace{-0.15cm}+(\mathbf{x}(k)-\mathbf{x}(k-1))^T (H \otimes I_N) (\mathbf{x}(k)-\mathbf{x}^*)\nonumber \\
&&\hspace{-0.15cm}+  \mathbf{x}(k\!-\!1\!)^T {c_0}  \mathbf{L} (\mathbf{x}(k)-\mathbf{x}^*) =0.\nonumber
\end{eqnarray}  
Combine the first {and fifth} term in \eqref{optim_equation_Nash_combo_multiplic_Revised_SumVec} and write,  \vspace{-0.35cm} 
\begin{eqnarray}\label{optim_equation_Nash_combo_multiplic_Revised_SumVec1}
&&\hspace{-0.15cm}\big(\!\mathbf{R} \mathbf{F}(\mathbf{x}(k\!-\!1\!)\!) \!+\! {c_0} \mathbf{L} \mathbf{x}(k\!-\!1\!)\!-\!\mathbf{R}\mathbf{F}(\mathbf{x}^*\!) \big)^T\!\!(\mathbf{x}(k)\!-\!\mathbf{x}^*)\\
&&\hspace{-0.15cm}=\!\big(\mathbf{R}\mathbf{F}(\mathbf{x}(k\!-\!1\!)\!)\! -\!\mathbf{R}\mathbf{F}(\mathbf{x}^*\!) +\! {c_0} \mathbf{L} \mathbf{x}(k\!-\!1\!)\!\big)^T\!\!(\mathbf{x}(k\!-\!1)\!-\!\mathbf{x}^*)\nonumber\\
&&\hspace{-0.15cm}+\big(\mathbf{R} \mathbf{F}(\mathbf{x}(k\!-\!1))\!-\! \mathbf{R} \mathbf{F}(\mathbf{x}^*\!)+\! {c_0} \mathbf{L} \mathbf{x}(k\!-\!1\!)\!\big)^T\!\!\!(\mathbf{x}(k)\!-\!\mathbf{x}(k\!-\!1)).\nonumber
\end{eqnarray}
{For the first term on the right-hand side of \eqref{optim_equation_Nash_combo_multiplic_Revised_SumVec1} we use  \eqref{eqboldRF_plus_cL_strgmon} in Lemma  \ref{boldRF_plus_cL_strgmon_lemma} (with $ \mathbf{L} \mathbf{x}^*=\zero_{N^2}$), 
while for the second term we use \eqref{inequality} for  $\rho=\frac{\bar{\theta}^2}{2\bar{\mu}}$,  
$\bar{\theta}=\theta \!+\! 2 c_0 d^*$. 
This yields, } \vspace{-0.3cm} 
\begin{eqnarray}\label{Cocoe_in_equatVec00}
&&\hspace{-0.15cm}\!\big(\mathbf{R} \!\mathbf{F}(\mathbf{x}(k\!-\!1\!))\!-\! \mathbf{R}\mathbf{F}(\mathbf{x}^*) \!+\! { c_0}  \mathbf{L} \mathbf{x}(k\!-\!1\!)\!\big)^T\! \!(\mathbf{x}(k)\!-\!\mathbf{x}^*)\nonumber\\
&&\hspace{-0.15cm}\geq{\bar{\mu}\|\mathbf{x}(k-1)-\mathbf{x}^*)\|^2} \nonumber \\
&&\hspace{-0.07cm}- {\frac{\bar{\mu}}{\bar{\theta}^2}} \| \mathbf{R}\mathbf{F}(\mathbf{x}(k\!-\!1))\!-\!\mathbf{R}\mathbf{F}(\mathbf{x}^*) \!+\! { c_0}  \mathbf{L} \mathbf{x}(k\!-\!1\!) \|^2\!
 \\
&&\hspace{-0.07cm}-\frac{\bar{\theta}^2}{4{\bar{\mu}} }\|\mathbf{x}(k)\!-\!\mathbf{x}(k\!-\!1)\|^2. \nonumber 
\end{eqnarray}
{Note that $ \mathbf{L}$ is $ \|\mathbf{L} \|_2$-Lipschitz with $\| \mathbf{L}\|_2 = \lambda_{max}(L) \leq 2 d^*$. Based on this, with $ \mathbf{L} \mathbf{x}^*=\zero_{N^2}$, {Assumption \ref{Lipchitz_boldFassump}} for $\mathbf{F}$, $\| \mathbf{R} \|_2 =1$ and the triangle inequality, we can write for the second term on the right-hand side of \eqref{Cocoe_in_equatVec00}, } \vspace{-0.2cm}
$$
 {\| \mathbf{R}\mathbf{F}(\mathbf{x}(k\!-\!1))\!-\!\mathbf{R}\mathbf{F}(\mathbf{x}^*) \!+\! { c_0}  \mathbf{L} \mathbf{x}(k\!-\!1\!) \|
 \! \leq \bar{\theta} \|\mathbf{x}(k-1)-\mathbf{x}^*)\|}
$$
{where $\bar{\theta}=\theta + 2 c_0 d^*$. Using this in  \eqref{Cocoe_in_equatVec00} yields} \vspace{-0.2cm}
\begin{eqnarray}\label{Cocoe_in_equatVec}
&&\hspace{-0.15cm}\!\big(\mathbf{R} \!\mathbf{F}(\mathbf{x}(k\!-\!1\!))\!-\! \mathbf{R}\mathbf{F}(\mathbf{x}^*) \!+\! { c_0}  \mathbf{L} \mathbf{x}(k\!-\!1\!)\!\big)^T\! \!(\mathbf{x}(k)\!-\!\mathbf{x}^*)\nonumber\\
&&\hspace{-0.15cm}\geq
 -\frac{{\bar{\theta}^2}}{4{\bar{\mu}} }\|\mathbf{x}(k)-\mathbf{x}(k-1)\|^2
\end{eqnarray}
Back to \eqref{optim_equation_Nash_combo_multiplic_Revised_SumVec}, for the second  term, since {$y \in   G(x(k))$, $y^* \in   G(x^*)$}, 
   $\mathcal{I}_{\Omega_i}$  is convex (Assumption~\ref{assump})  it follows that  \vspace{-0.3cm} 
 \begin{equation}\label{convexity_of_IndicatorVec}
(y-y^*)^T(x(k)-x^*)\geq 0. 
\end{equation}
Finally for the third  term  in \eqref{optim_equation_Nash_combo_multiplic_Revised_SumVec}, we use the following. 
From 
{\eqref{wstep_Algvector} with zero initial conditions, }
it follows that $\mathbf{w}(k+1)= c \mathbf{L} \mathbf{q}(k)$,  where $\mathbf{q}(k) = \sum_{t=0}^k \mathbf{x}(t)$.
Thus,   \vspace{-0.3cm} 
\begin{equation}\label{delta_q}
\mathbf{q}(k) - \mathbf{q}(k-1) =  \mathbf{x}(k)  
\end{equation}
Then, we can write  \vspace{-0.3cm} 
\begin{eqnarray*}\label{simpil_u_underline_befVecAlt}
&&\hspace{-0.15cm} (\mathbf{w}(k+1)-{c \mathbf{L}  \mathbf{q}^*})^T(\mathbf{x}(k)-\mathbf{x}^*) \\
&&\hspace{-0.15cm}  =(\mathbf{q}(k)-  \mathbf{q}^*)^T c\mathbf{L}(\mathbf{q}(k) - \mathbf{q}(k-1)-\mathbf{x}^*) \nonumber \\ 
&&\hspace{-0.15cm}  = (\mathbf{q}(k)- \mathbf{q}^*)^Tc\mathbf{L}(\mathbf{q}(k) - \mathbf{q}(k-1)) \nonumber
\end{eqnarray*}
Using  \eqref{Cocoe_in_equatVec}, \eqref{convexity_of_IndicatorVec} and
  this last relation  in \eqref{optim_equation_Nash_combo_multiplic_Revised_SumVec} yields\vspace{-0.25cm}
\begin{eqnarray}\label{optim_equation_Nash_combo_multiplic_Revised_SumVecINEQ1Alt}
&&\hspace{-0.15cm}\!-\frac{{\bar{\theta}^2}}{4{\bar{\mu}}}\|\mathbf{x}(k) \!- \!\mathbf{x}(k \!-\!1)\|^2  
\!+\!(\mathbf{q}(k) \!- \!  \mathbf{q}^*\!)^T c \mathbf{L} (\!\mathbf{q}(k) \!-\! \mathbf{q}(k\!-\!1)\!) \nonumber \\
&&\hspace{-0.0cm}+(\mathbf{x}(k)-\mathbf{x}(k-1))^T (H \otimes I_N) (\mathbf{x}(k)-\mathbf{x}^*) \leq 0,
\end{eqnarray}
where  $H\otimes I_N{\succ}0$. 
Using \eqref{simple_math} to deal with all cross-terms in the previous inequality, yields \vspace{-0.25cm}
\begin{eqnarray*}\label{optim_equation_Nash_combo_multiplic_Revised_Sum_simpil_Q>0Alt}
&&\hspace{-0.15cm} \! \| \mathbf{x}(k) \!-\! \mathbf{x}^*\|^2_{H \!\otimes  \!I_N} \! \!+ \!  \| \!\mathbf{q}(k) \!- \! \mathbf{q}^*\|^2_{c \mathbf{L}} 
\!-\!  \| \mathbf{x}(k\!-\!1) \!- \! \mathbf{x}^*\|^2_{H \! \otimes \! I_N} \! 
\nonumber  \\
&&\hspace{1cm}  - \!  \| \mathbf{q}(k\!-\!1) \!-\! \mathbf{q}^*\|^2_{c\mathbf{L}} 
\, \leq 
\, - \! \| \mathbf{x}(k) \!-\! \mathbf{x}(k\!-\!1) \|^2_{H \!\otimes \! I_N}  \! \!
\nonumber  \\
&&\hspace{1.5cm} 
 - \!  \| \mathbf{q}(k) \!- \! \mathbf{q}(k \!- \!1) \|^2_{c\mathbf{L}} \!
 \!+ \! \frac{{\bar{\theta}^2}}{2{\bar{\mu}}} \! \| \!\mathbf{x}(k) \!- \!\mathbf{x}(k \!-\!1) \!\|^2.   \nonumber
\end{eqnarray*}
Let us make  the following notations:  \vspace{-0.3cm} 
\begin{equation}\label{tafkik_z}
\bz(k)\!:=\!\left [ \!
\begin{array}{c}
\bx(k) \\
\bq(k) \\
\end{array} \!
\right ]\!, \, 
\bz^* \!:= \!\left [ \!
\begin{array}{c}
\mathbf{x}^* \\  
\mathbf{q}^* \\
\end{array} \!
\right ]\!,\CG \!=\! \left [ \!   
\begin{array}{cc}
H \hspace{0.5cm}&  \mathbf{0}_{N \times N} \\
\mathbf{0}_{N \times N} & cL \\
\end{array} \!
\right ]
\end{equation}
Then we can write the last inequality as   \vspace{-0.3cm}  
\begin{eqnarray*}\label{optim_equation_Nash_combo_multiplic_Revised_Sum_simpil_Q>0AltBis}
&&\hspace{-0.15cm}  \| \mathbf{z}(k)- \mathbf{z}^*\|^2_{\CG \otimes I_N} -  \| \mathbf{z}(k-1)- \mathbf{z}^*\|^2_{\CG \otimes I_N}   \\
&&\hspace{-0.15cm} \leq - \| \mathbf{z}(k)- \mathbf{z}(k-1) \|^2_{\CG\otimes I_N}   + \frac{{\bar{\theta}^2}}{2{\bar{\mu}}}\|\mathbf{x}(k)-\mathbf{x}(k-1)\|^2.   \nonumber
\end{eqnarray*}
Moreover, since $ \frac{{\bar{\theta}^2}}{2{\bar{\mu}}} \|\mathbf{x}(k)-\mathbf{x}(k-1)\|^2  \leq  \frac{{\bar{\theta}^2}}{2{\bar{\mu}}\lambda_{\min}(H)} \|\mathbf{x}(k)-\mathbf{x}(k-1)\|^2_{H \otimes I_N}$
$ \leq \frac{{\bar{\theta}^2}}{2{\bar{\mu}} \lambda_{\min}(H)}\|\mathbf{z}(k)-\mathbf{z}(k-1)\|^2_{\CG \otimes I_N} $, 
from  the foregoing it follows that   \vspace{-0.25cm}
\begin{eqnarray}\label{optim_equation_Nash_combo_multiplic_Revised_Sum_simpil_Q>0AltBisBis}
&&\hspace{-0.15cm}  \| \mathbf{z}(k)- \mathbf{z}^*\|^2_{\CG \otimes I_N} -  \| \mathbf{z}(k-1)- \mathbf{z}^*\|^2_{\CG \otimes I_N}  \leq \\
&&\hspace{2cm} - \zeta \| \mathbf{z}(k)- \mathbf{z}(k-1) \|^2_{\CG\otimes I_N}    \leq 0,   \nonumber
\end{eqnarray}
where $\zeta = 1 -  \frac{{\bar{\theta}^2}}{2{\bar{\mu}} \lambda_{\min}(H)},$ and $0 < \zeta <1$ by \eqref{condition}. Summing \eqref{optim_equation_Nash_combo_multiplic_Revised_Sum_simpil_Q>0AltBisBis} over $k$ from 1 to $\infty$ yields  \vspace{-0.35cm} 
\begin{eqnarray}\label{optim_equation_Nash_combo_multiplic_Revised_Sum_simpil_Q>0AltSumk}
&&\hspace{-0.15cm}\sum_{k=1}^\infty \|\mathbf{z}(k)-\mathbf{z}(k-1)\|_{\CG \otimes I_N}^2 
\hspace{-0.1cm}\leq\hspace{-0.1cm} \frac{1}{\zeta} \|\mathbf{z}(0)-\mathbf{z}^*\|_{\CG \otimes I_N}^2 \hspace{-0.1cm} < \hspace{-0.05cm}\infty  
\end{eqnarray}
From \eqref{optim_equation_Nash_combo_multiplic_Revised_Sum_simpil_Q>0AltSumk} it follows that $ \|\mathbf{z}(k)-\mathbf{z}(k-1)\|^2_{\CG\otimes I_N}\rightarrow 0$. 
From \eqref{optim_equation_Nash_combo_multiplic_Revised_Sum_simpil_Q>0AltBisBis}, since $\|\mathbf{z}(k)-\mathbf{z}^*\|_{\CG\otimes I_N}^2$   
{ or  $ \| \mathbf{x}(k)- \mathbf{x}^*\|^2_{H \otimes I_N} +  \| \mathbf{q}(k)- \mathbf{q}^*\|^2_{c\mathbf{L}}$ is bounded and non-increasing, it follows that 
$ \| \mathbf{x}(k)- \mathbf{x}^*\|^2_{H \otimes I_N} +  \| \mathbf{q}(k)- \mathbf{q}^*\|^2_{c\mathbf{L}} \rightarrow  v$, for some  $v \geq0$, 
where  $(H \otimes I_N) \succ 0$,   $ \mathbf{L} \succeq 0$. }

\vspace{-0.3cm}
{ Any $\mathbf{q} \in \mathbb{R}^{N^2}$ can be decomposed as $\mathbf{q} = \mathbf{q}^{\|} + \mathbf{q}^{\perp}$, where $ \mathbf{q}^{\|} \in Ker (\mathbf{L})$, 
$ \mathbf{q}^{\perp} \in Ker(\mathbf{L})^{\perp}$,  
with  $({\mathbf{q}^{\perp}})^T \mathbf{L}  \mathbf{q}^{\perp} >0$, $\forall \mathbf{q}^{\perp} \neq 0$. 
Decomposing  $\mathbf{q}(k)$ and  $\mathbf{q}^*$ into $(\cdot)^{\|}$ and $(\cdot)^{\perp}$ components,  $  \| \mathbf{q}(k)- \mathbf{q}^*\|^2_{c\mathbf{L}} =  \| \mathbf{q}(k)^{\perp}- \mathbf{q}^{*\perp}\|^2_{c\mathbf{L}}$, with $\mathbf{q}(k)^{\perp}, \mathbf{q}^{*\perp} \in Ker(\mathbf{L})^{\perp}$. 
Hence, from the above, 
\vspace{-0.25cm}
$$V_{x^*}(k):= \| \mathbf{x}(k)- \mathbf{x}^*\|^2_{H \otimes I_N} +  \| \mathbf{q}(k)^{\perp}- \mathbf{q}^{*\perp}\|^2_{c\mathbf{L}}$$ is bounded and non-increasing and converges to $v$.  Thus, since  $H \succ 0$ and  $ \| \mathbf{q}(k)^{\perp}- \mathbf{q}^{*\perp}\|^2_{c\mathbf{L}} >0 $,  $\forall \mathbf{q}^{\perp}(k) \neq 0,$   
the sequence $[\mathbf{x}(k) \,; \, \mathbf{q}(k)^{\perp}]$ is bounded, hence has at least a limit  point 
$[\mathbf{\overline{x}} \; ; \, \mathbf{\overline{q}}^{\perp}]$, where 
$\mathbf{\overline{q}}^{\perp} \in Ker(\mathbf{L})^{\perp}$.  
{Also, $\mathbf{w}(k+1)= c\mathbf{L} \mathbf{q}(k) = c\mathbf{L} \mathbf{q}(k)^{\perp}$, 
has a limit point $\mathbf{\overline{w}} = c \mathbf{L} \mathbf{\overline{q}}^{\perp}$. }  
By Lemma~\ref{correctness_lemma_1},  any limit point  
of $\mathbf{x}(k)$ satisfies $\mathbf{\overline{x}} = \mathbf{1}_N \otimes x^*$, with $x^*$ 
  NE of the game  and $\mathbf{\overline{w}} $ such that \eqref{combo_no_stackBFG} holds, 
  hence $\mathbf{\overline{x}} = \mathbf{x}^*$. 
Consider  now\vspace{-0.25cm}  
\begin{equation}\label{eq:lastV}
V_{\overline{x}}(k):= \| \mathbf{x}(k)- \mathbf{x}^*\|^2_{H \otimes I_N} +  \| \mathbf{q}(k)^{\perp}- \overline{\mathbf{q}}^{\perp}\|^2_{c \mathbf{L}}
\end{equation}
Similar to the first part of the proof for   $V_{x^*}(k)$ with $ \mathbf{x}^*$, $ c \mathbf{L} \mathbf{q}^*$  { and \eqref{combo_no_stackBFGVec_ZERO_q}}, it can be shown that $V_{\overline{x}}(k)$  
 is bounded and non-increasing, by using   \eqref{combo_no_stackBFG} for $\mathbf{\overline{x}} = \mathbf{x}^*$, $\mathbf{\overline{w}} = c \mathbf{L} \mathbf{\overline{q}}^{\perp}$ instead of \eqref{combo_no_stackBFGVec_ZERO_q}.  
Hence, $V_{\overline{x}}(k)$ converges to some $\overline{v} \geq 0$ as $k \rightarrow \infty$. This $\overline{v}$ is the same along any subsequence. Since $[\mathbf{x}^* \; ; \, \mathbf{\overline{q}}^{\perp}]$ is a limit point for $[\mathbf{x}(k) \,; \, \mathbf{q}(k)^{\perp}]$,  there exists a subsequence $\{  k_n \}$,  $ k_n  \rightarrow \infty$  as $n \rightarrow \infty$, such that  $\{ [\mathbf{x}(k_n) \,; \, \mathbf{q}(k_n)^{\perp}] \} \rightarrow [\mathbf{x}^* \; ; \, \mathbf{\overline{q}}^{\perp}]$ as $n \rightarrow \infty$. 
 Taking limit  of  $V_{\overline{x}}(k_n)$, \eqref{eq:lastV}, along  this $\{ k_n\}$ subsequence, as $n \rightarrow \infty$, it follows that $\overline{v} =0$. 
Since  $(H \!\otimes  \! I_N) \!\succ 0$ and $\mathbf{\overline{q}}^{\perp}, \mathbf{q}(k)^{\perp} \in Ker(\mathbf{L})^{\perp}$,  taking limit in  \eqref{eq:lastV} as $k \rightarrow \infty$,  it follows that $\mathbf{x}(k) \!\rightarrow \mathbf{x}^*$,  with $\mathbf{x}^*\! =\! \mathbf{1}_N \!\otimes \!x^*$ and $x^*$ NE. 
$\hfill\blacksquare$

\begin{remark}\label{vsAuto2007}
We note that, unlike ADMM for DOP, \cite{WeiOzdaglar2013arX}, \cite{HongLuoRaza2016}, to prove convergence we cannot use a common objective function and common Lagrangian. Instead our proof relies on properties of the extended pseudo-gradient  $\mathbf{F}$. {Compared  to \cite{salehiPavelIFAC2017}, \cite{WeiShiACC}, where $c_0\!=\!0$ and cocoercivity of $\mathbf{R}\mathbf{F}$ was assumed, here we only assume Lipschitz continuity  of  $\mathbf{F}$. 
 We achieved this relaxation by slightly modifying the ADMM algorithm so that the $\mathbf{x}$-update uses the extra $c_0\!>\!0$ parameter in the proximal-approximation (quadratic penalty). Convergence is proved for sufficiently large $c_0$.  This is  similar to the way in which non-convexity is overcome in ADMM for DOP  in \cite{HongLuoRaza2016}, \cite{WangADMM2015}). 
Based on Gershgorin theorem, a sufficient diagonal-dominance condition can be derived for \eqref{condition} to hold, and   $\beta_i$ parameters  can be selected independently by players (see e.g. \cite{yipengLP1}).  We note that assumptions on $F$ and $\mathbf{F}$ could be relaxed to hold only locally around $x^*$  and $\mathbf{x}^*$, in which case all results become local. 
  We also note that the class of quadratic games satisfies all asumptions globally.
 }
\end{remark}

\vspace{-0.3cm}
\subsection{Convergence Rate Analysis} 


 \vspace{-0.3cm}
Next we investigate the convergence rate of  Algorithm~\ref{ADMMalgorithm}. We use the following result. \vspace{-0.25cm}
\begin{proposition}[\cite{shi2015proximal}]\label{prop:O_1_k}
	If a sequence $\{a(k)\}\subset\R$ is: \emph{(i)} nonnegative, $a(k)\geq0$, \emph{(ii)} summable, $\sum_{k=1}^{\infty}a(k)<\infty$, and \emph{(iii)} monotonically non-increasing, $a({k+1})\leq a(k)$, then we have: $a(k)=o\left(\frac{1}{k}\right)$, i.e., $\lim_{k\rightarrow\infty} k \, a(k)=0$.
\end{proposition} \vspace{-0.3cm}
For  $\|\bz(k-1)-\bz(k)\|_{\CG \otimes I_N}^2$,  with $\bz(k)$, ${\CG \otimes I_N}$ as in \eqref{tafkik_z}, \eqref{tafkik_z},  summability follows from \eqref{optim_equation_Nash_combo_multiplic_Revised_Sum_simpil_Q>0AltSumk} in Theorem~\ref{theorem_convergence_Algo1}. Next   we provide a lemma showing that $\|\bz({k-1})-\bz(k)\|_{\CG \otimes I_N}^2$ is monotonically non-increasing.  
\begin{lemma}\label{lemma:monotonicity} 
	Under the  assumptions of Theorem~\ref{theorem_convergence_Algo1}, 
	 the sequence $\{\bz(k)\}$, \eqref{tafkik_z}, generated by Algorithm~\ref{ADMMalgorithm} 
	 satisfies for all $ k\geq1$, \vspace{-0.3cm}
	\begin{equation*}\label{eq:monotonicity_PG}
	\hspace{0em}\begin{array}{rcl}
	\|\bz(k)-\bz(k+1)\|_{\CG \otimes I_N}^2\leq\|\bz(k-1)-\bz(k)\|_{\CG \otimes I_N}^2,
	\end{array} \vspace{-0.25cm}
	\end{equation*}
	where ${\CG \otimes I_N}\succeq0$ is as in \eqref{tafkik_z}.
\end{lemma} \vspace{-0.25cm}
\par{\emph{Proof}}. 
Let $\D\bx(k+1):=\bx(k)-\bx(k+1)$, $\D\bq(k+1):=\bq(k)-\bq(k+1)$, $\D\bz(k+1):=\bz(k)-\bz(k+1)$, $\D\mathbf{F}(\bx(k+1)):=\mathbf{F}(\bx(k))-\mathbf{F}(\bx(k+1))$, $\D\mathbf{G}(x(k+1))\:=\mathbf{G}(x(k))-\mathbf{G}(x(k+1))$. 
Following similar arguments as in the proof of  Theorem~\ref{theorem_convergence_Algo1}, we obtain  similar to \eqref{optim_equation_Nash_combo_multiplic_Revised_Sum_simpil_Q>0AltBisBis}, 
$\|\D\bz(k)\|_{\CG \otimes I_N}^2 \!-\!\|\D\bz({k+1})\|_{\CG \otimes I_N}^2
\! \geq \! \zeta\|\D\bz(k)\!-\!\D\bz(k+1)\|_{\CG \otimes I_N}^2 \!\geq \!0, 
$ 
where $\zeta=1-\frac{{\bar{\theta}^2}}{2{\bar{\mu}}\lambda_{\min}(H)}>0,$ by \eqref{condition}. 
$\hfill\blacksquare$ \vspace{-0.3cm}

\begin{theorem}\label{theorem:o_1_k} 
	Under the same assumptions of Theorem~\ref{theorem_convergence_Algo1}, the
	following rate{s} hold for  Algorithm~\ref{ADMMalgorithm}:  \vspace{-0.3cm} 
	\begin{equation*}\label{eq:suc_dif_1_k_x}
	\begin{array}{l}
	\|\bx{^\perp}(k)\|_{\mathbf{L}}^2=o\left(\frac{1}{k}\right), 
	\end{array} \vspace{-0.3cm}
	\end{equation*}
	where {$\bx^\perp(k) \in Ker(\mathbf{L})^{\perp}$}, $\mathbf{L}=L \otimes I_N$,  $L$ the Laplacian matrix of $G_c$, and  for some $y \in \mathbf{G}(x(k))$,	
 \vspace{-0.3cm}
	$$\hspace{-0.25cm}
\!\| \mathbf{R} (\!\mathbf{\mathbf{F}}( \!\mathbf{x}(k\!-\!1\!))\!+  y \!)\!+\! \mathbf{w}\!(k\!+\!1\!) 
\!+\! c_0 \mathbf{L} \mathbf{x}^\perp(k\!-\!1\!) \|^2_{\!(H \otimes I_N)^{\!-1}}\!=\! o\left(\!\frac{1}{k}\! \right)  
$$ \vspace{-0.4cm}
\end{theorem}\vspace{-0.3cm}
\par{\emph{Proof}}.
	By Theorem~\ref{theorem_convergence_Algo1} (equation \eqref{optim_equation_Nash_combo_multiplic_Revised_Sum_simpil_Q>0AltSumk}), 
	$\sum_{k=0}^\infty\|\bz(k-1)-\bz(k)\|_{\CG \otimes I_N}^2<\infty$. Moreover, Lemma~\ref{lemma:monotonicity} proves that $\|\bz(k-1)-\bz(k)\|_{\CG \otimes I_N}^2$ is monotonically non-increasing. Then, it directly follows by Proposition \ref{prop:O_1_k} that, \vspace{-0.3cm}
	\begin{equation*}\label{eq:suc_dif_1_k_z}
	\begin{array}{l}
	\|\bz(k-1)-\bz(k)\|_{\CG \otimes I_N}^2=o\left(\frac{1}{k}\right).
	\end{array} \vspace{-0.3cm}
	\end{equation*}
Expanding this by using \eqref{tafkik_z} yields, \vspace{-0.3cm}
$$\|\bx(k-1)-\bx(k)\|_{H \otimes I_N}^2+c\|\bq(k-1)-\bq(k)\|_{\mathbf{L}}^2 =o\left(\frac{1}{k}\right) \vspace{-0.3cm}
$$
and, using \eqref{delta_q} $\bq(k-1)-\bq(k) = \bx(k)$, so that  \vspace{-0.3cm}
$$\|\bx(k-1)-\bx(k)\|_{H \otimes I_N}^2+c\|\bx(k)\|_{\mathbf{L}}^2 =o\left(\frac{1}{k}\right) \vspace{-0.3cm}
$$
As in the proof of Theorem~\ref{theorem_convergence_Algo1}, we can write this as \vspace{-0.3cm} 
\begin{equation}\label{esbat_convergence_rate_optimality1}
\|\bx(k-1)-\bx(k)\|_{H \otimes I_N}^2+c\| \bx^\perp(k)\|_{\mathbf{L}}^2 =o\left(\frac{1}{k}\right) 
\end{equation}
where  $\bx^\perp(k) \in Ker(\mathbf{L})^{\perp}$, $\bx(k)^\perp \mathbf{L} \bx^\perp(k) >0$, for all $ \bx^\perp(k) \neq 0$.   
 If the sum of two positive sequences is $o(\frac{1}{k})$, then each sequence is $o(\frac{1}{k})$, hence, by  the foregoing, 
 \vspace{-0.25cm}
\begin{equation*}\label{esbat_convergence_rate_optimality_tahesh}
\|\bx{^\perp}(k)\|_{\mathbf{L}}^2=o\left(\frac{1}{k}\right).
\end{equation*} 
From \eqref{esbat_convergence_rate_optimality1}, \eqref{step67_altVec}, we can write 
$\!
 \| \mathbf{x}(k)\!-\!  \mathbf{x}(k\!-\!1) \|^2_{H \otimes I_N} 
\!=\! \| \mathbf{R} \!\mathbf{\mathbf{F}}( \mathbf{x}(k\!-\!1))\!+\! \mathbf{R} y \!+\! \mathbf{w}(k\!+\!1)\! +\!{ c_0} \mathbf{L} \mathbf{x}^\perp(k\!-\!1\!) \|^2_{(H \otimes I_N)^{-1}}\!\!=o\left(\frac{1}{k}\right)
$, 
for some $y \in \mathbf{G}(x(k))$. 
$\hfill\blacksquare$

 \vspace{-0.3cm}
 {\begin{remark}
Since $\bx^\perp(k) \in Ker(\mathbf{L})^{\perp}$,  
$ \| \bx^\perp(k)\|_{\mathbf{L}}^2$ can be interpreted as a weighted distance from the consensus subspace, while  
the second rate  can be interpreted a weighted distance from the optimality conditions. 
Unlike DOP, e.g. \cite{WeiOzdaglar2013arX}, in a game we do not have a common merit function for all players and the convergence rate we are able to show is only in terms of optimality residuals. 
  \vspace{-0.3cm}
\end{remark}
}
\section{Simulation Results}\label{sec:simulations}
In this section, we  present two examples to illustrate the performance of Algorithm~\ref{ADMMalgorithm}. The first one is a networked Nash-Cournot game (Example \ref{Ex_quadratic}) as in \cite{Jayash}, while the second one is a rate control game in Wireless Ad-Hoc Network (WANET),  as  in \cite{salehisadaghianiPavel2018Auto,Alpcan2005}. \vspace{-0.3cm}

\subsection{Example 1: Nash-Cournot game over networks}\vspace{-0.3cm}



In this section, we examine the performance of the proposed algorithm on a networked Nash-Cournot game (Example \ref{Ex_quadratic}).  We consider $N=20$ firms and $m=7$ markets, 
(Fig. \ref{fig_network_cournot_game}). 
Each firm $i$ has a local normalized constraint as $0\leq x_{T,i} \leq 1$ on its total production $x_{i,T}=n_i x_i$,  
where $n_i$ is the number of markets it participates in, hence $\Omega_i = [0,1/n_i]$. The local objective function is taken as \eqref{network_cournot_game_function}, where the local production cost function of firm $i$ is $c_i(x_i)=  {n_i}^2 q_i x_{i}^2+ n_i  b_i x_i$ where  $q_i$, $b_i$ are randomly drawn from  the uniform distribution $U(1,2)$. All  assumptions hold globally. In the first case,  each local cost function $J_i$ is normalized by dividing by $n_i^2 q_i$, which yields a better condition number. The price function is taken as the linear function $ P=\bar{P}-Z Ax$, where  $z_k =0.01$ and  $\bar{P}_k$  are randomly drawn from $ZA\Omega_i+U(1,2)$.  This yields $\lambda_{\min}(Q) = 1.001=\mu$,  $\lambda_{\max}(Q) = 1.099=\theta$. For a communication graph  $G_c$ as in Fig. \ref{fig_G_c_graph}, $\lambda_2(L) = 0.102$, and the bound for $c_0$ in Lemma \ref{boldRF_plus_cL_strgmon_lemma} is $22.5$. We set   $c_0 = 22.5$, an arbitrary step-size $c=1$, and $\beta_i = 10 $ for all firms, satisfying \eqref{condition}. 
The results for the implementation of Algorithm~\ref{ADMMalgorithm}  in the partial decision setting, where the primal variables $x^i$ are exchanged over  the sparsely connected graph $G_c$ (Fig.~\ref{fig_G_c_graph}),  are shown in Fig. \ref{eps:graph_and_multi_gossip}. Initial conditions are randomly selected, for $x_i$ from $U(0,0.5)$, and for  $x_j^i$ ($i\neq j$) from $U(0,1)$. The results indicate fast convergence.   

\begin{figure}
		\vspace{-3.0cm}
		\includegraphics [scale=0.4]{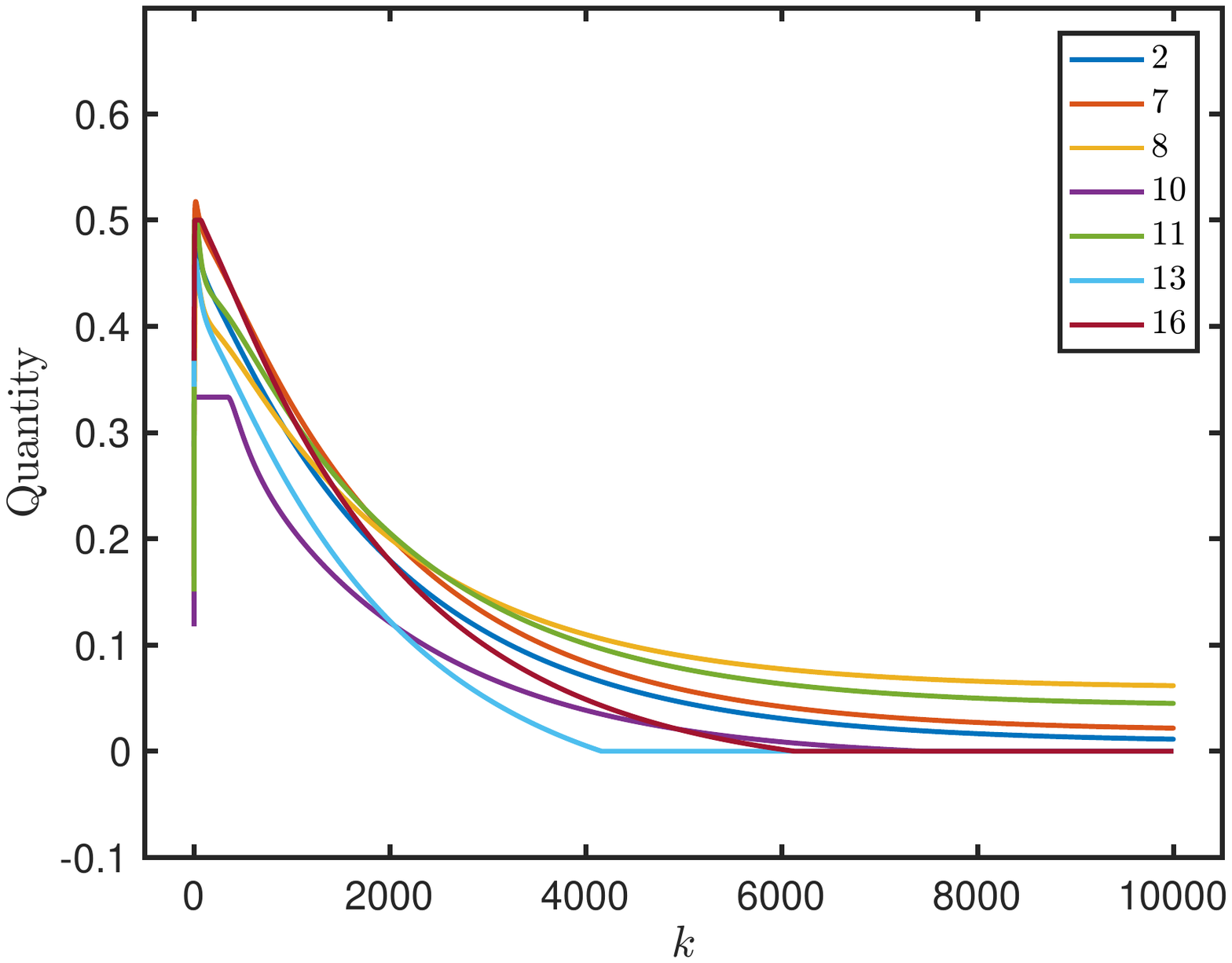}
		\vspace{-3.6cm}
	\caption{NE computation of a Nash-Cournot game by Algorithm~\ref{ADMMalgorithm} over $G_c$ in Fig.~\ref{fig_G_c_graph}. 
	The figure shows the convergence paths of 7 randomly selected firms.}\label{eps:graph_and_multi_gossip}\vspace{-0cm}
\end{figure}\vspace{-0.25cm}

Next,  we study an extreme case where the costs are widely different. This yields a high condition number, and we can study how the cost functions affect the outcome of market competition. The  $q_i$ and $b_i$ for the first $10$ firms are set at the level of $0.01$ while those for the other firms are at the level of $100$, and each firm $i$ has a local constraint $0 \leq x_{i,T} \leq 100$.  Intuitively, this will lead to high/low market occupation of the first/last $10$ firms. The price function is taken as the linear function $ P=\bar{P}-ZAx$, where $z_k$ are randomly drawn from  $U(1,10)$ and $\bar{P}_k$ randomly generated such that $P_k>0$ for all feasible $x$. 
 This yields $\lambda_{\min}(Q) = 1.482=\mu$,  $\lambda_{\max}(Q) = 1005=\theta$, and 
  the bound for $c_0$ in Lemma \ref{boldRF_plus_cL_strgmon_lemma} is $6e6$. However, we found this bound conservative in our experiments. We plot the sequences of relative error (vs. the 'true' NE obtained by (centralized) gradient method) under different parameter settings in the top sub-figure of Fig. \ref{eps:graph_and_multi_gossip_I}. 
  The trends reflected in the top sub-figure of Fig. \ref{eps:graph_and_multi_gossip_I} match the intuition of our algorithm design (see Appendix A) and the implication of \eqref{condition} in Theorem \ref{theorem_convergence_Algo1}. In general, smaller penalty parameters $c$, $c_0$, and $\beta_i$'s lead to faster convergence but they have to be large enough to guarantee convergence. An instance of the first $10^4$ iterations of trajectories of selected firms are plotted in the bottom sub-figure of Fig. \ref{eps:graph_and_multi_gossip_I}. The initial estimate for $x^i_{-i}$'s are randomly drawn from $(0,200)$; to observe diverse firm trajectories, the first $5$ and last $5$ firms' local $x^i_i$'s are initialized as the upper limits of their local constraints, respectively, while the other firms start from $0$. 
  This way, we can see four cases: (1-5) ``aggressive start, rich end"; (6-10) ``conservative start, rich end"; (11-15) ``conservative start, poor end" (16-20) ``aggressive start, poor end". 
  } 
\vspace{-0.2cm}
\begin{figure}
		\vspace{-3.0cm}
		\includegraphics [scale=0.4]{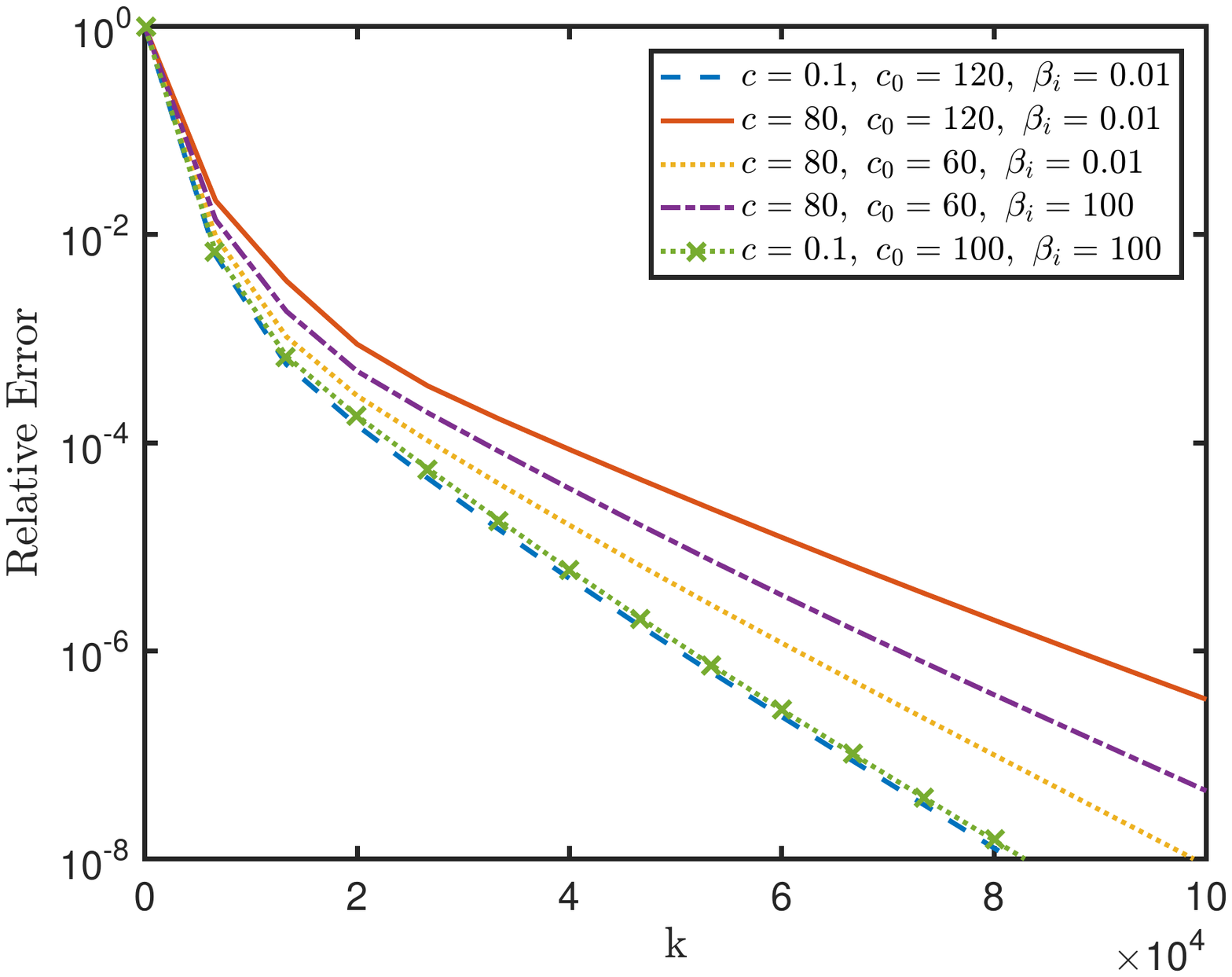}\vspace{-5.5cm}
		\includegraphics [scale=0.4]{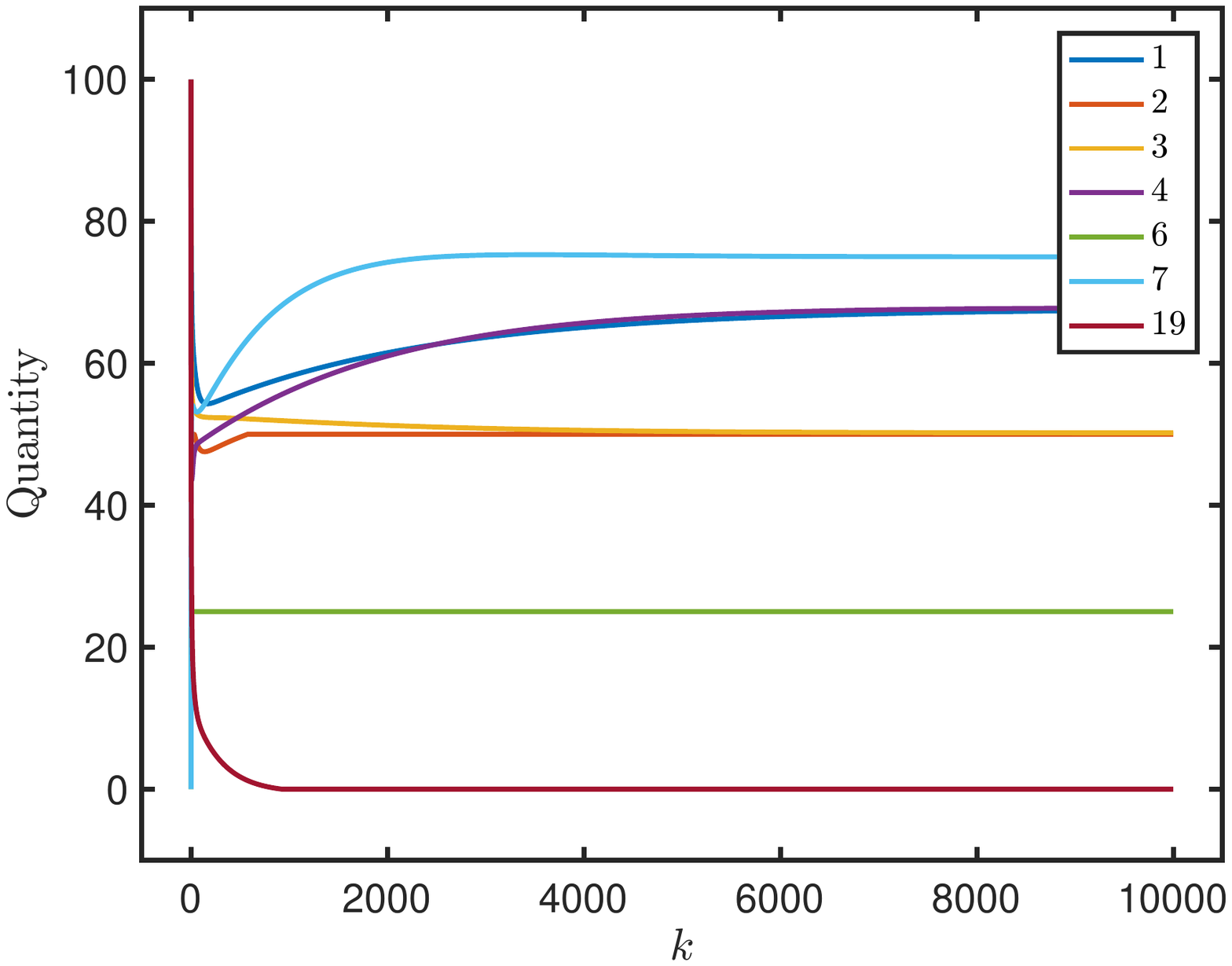}
		\vspace{-3.6cm}
		\caption{NE computation of a Nash-Cournot game by Algorithm~\ref{ADMMalgorithm} over $G_c$ in Fig.~\ref{fig_G_c_graph}. The top figure plots relative error $\frac{\|x(k)-x^*\|}{\|x(0)-x^*\|}$ vs. iteration number $k$. The bottom figure shows the convergence paths of 7 selected firms.
		}\label{eps:graph_and_multi_gossip_I}\vspace{-0cm}
\end{figure}\vspace{-0.15cm}


\subsection{Example 2: Rate control game over WANET}
\begin{figure}
 	\vspace{-2.25cm}
	\hspace{-2.23cm}
	\centering
	\includegraphics [scale=0.54]{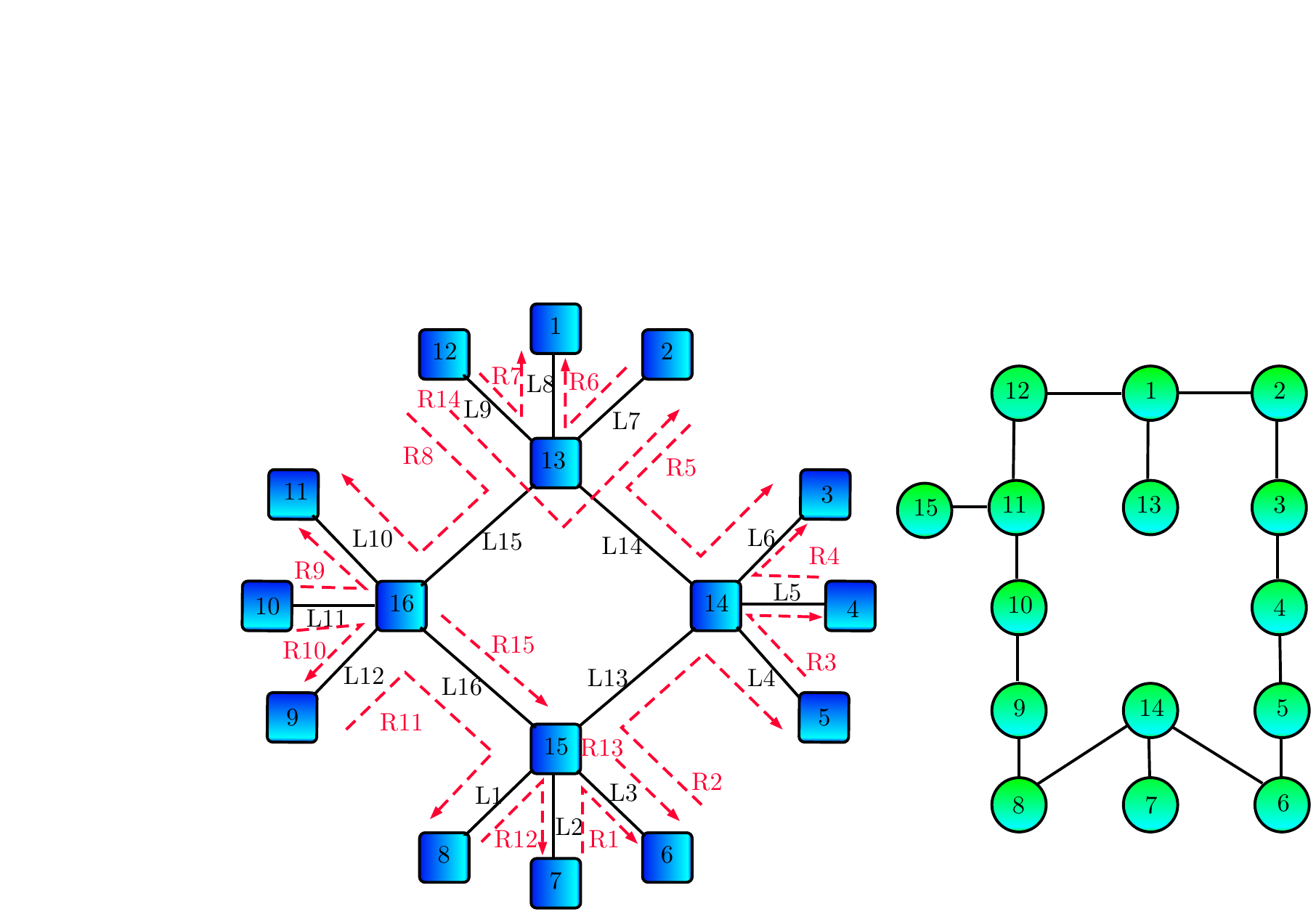}
	\vspace{-0.25cm}
	\caption{(a)  Wireless Ad-Hoc Network (left). (b) Communication graph $G_c$ (right) for Example 2.}
	\label{fig:minipage_a_b_FRED}	
\end{figure}

\begin{figure}
	\vspace{-3.58cm}
	\hspace{-0.8cm}
	\includegraphics [scale=0.45]{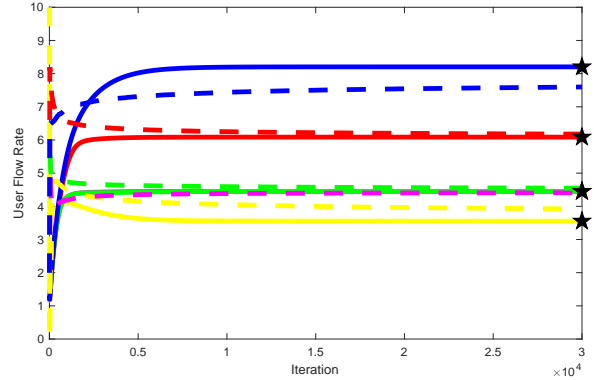}
	\vspace{-4.cm}
	\caption{Flow rates of users 1, 3, 5, 8, 13 using Algorithm~\ref{ADMMalgorithm} (solid lines) vs. algorithm in \cite{salehisadaghianiPavel2018Auto} (dashed lines). NE represented by black stars.}
\label{fig:minipage1_FRED}	
\end{figure}

\vspace{-0.2cm}
In this section we consider an engineering network-protocol example for a rate control game over  
a \emph{Wireless Ad-Hoc Network} (WANET), \cite{jayash7,Alpcan2005}. Consider  16 mobile nodes interconnected by multi-hop communication paths in a WANET Fig.~\ref{fig:minipage_a_b_FRED}~(a).  $N=15$ users/players want to transfer data from a source node to a destination node via this WANET.  In Fig.~\ref{fig:minipage_a_b_FRED}~(a) solid lines represent physical links and dashed lines display unique paths that are assigned to each user to transfer its data.  Each link has a positive capacity that restricts the users' data flow. {  In this WANET scenario there is no central node. Each user/player is prescribed to communicate locally with its neighbours over a preassigned (sparse but connected) communication graph $G_c$, e.g. as depicted in Fig.~\ref{fig:minipage_a_b_FRED}~(b). In Fig.~\ref{fig:minipage_a_b_FRED}~(b) a node represents a player,  and an edge in the graph signifies the ability for the two players to exchange information. Each user exchanges information  locally (thus mitigating the lack of centralized information), so that he  can estimate the other flows and, based on them, adjust its own decision which is his data rate.}
Here is the list of WANET notations:\vspace{-0.1cm}
\begin{enumerate}
	\item $L_j$: Link $j$, $j =1,\ldots,16$,
	\item $R_i$: The path assigned to user $i$, $i=1,\ldots,15$,
	\item $C_j>0$: Link   $j$'s capacity, $j =1,\ldots,16$, 
	\item $0\leq x_i\leq 10$: The data flow of user $i$, $i =1,\ldots,15$.
\end{enumerate}
Each path consists of a set of links, e.g., $R_1=\{L_2,L_3\}$.  Each user $i$ decision (action) is its data rate flow $0\leq x_i\leq 10$ to send over $R_i$, based on its cost function $J_i$, \vspace{-0.250cm}
\begin{equation*}\label{Cost_fcn_gen}
J_i(x_i,x_{-i}):=\sum_{j:L_j\in R_i}\frac{\kappa}{C_j-\sum_{w:L_j\in R_w}x_w}-\chi_i \log(x_i+1),\end{equation*}
which depends on the flows of the other users, and where $\kappa>0$ and $\chi_i>0$ are network-wide known and user-specific parameters, respectively. In this case  $\lambda_2(L) = 0.18$ and the other assumptions hold locally. For $\chi_i=10$ $\forall i=1,\ldots,15$ and $C_j=10$ $\forall j=1,\ldots,16$ we compute numerically $\|DF(x)\|_2$ and  $\|D\mathbf{F}(x)\|_2$  over $\Omega$,  $\theta_0= 4.7$, $\theta= 6.7$. The bound for $c_0$  in Lemma \ref{boldRF_plus_cL_strgmon_lemma} is $30.7$. We set   $c_0 = 31$,  hence $\bar{\mu} =0.04$, an arbitrary step-size $c=1$ and $\beta_i = 14$, $\forall  i=1,\ldots,15$ (by \eqref{condition} in Theorem \ref{theorem_convergence_Algo1}).  
We compare our Algorithm~\ref{ADMMalgorithm} with the gossip-based one proposed in \cite{salehisadaghianiPavel2018Auto} run  over $G_c$, with diminishing-step sizes. The results are shown in Fig.~\ref{fig:minipage1_FRED}, for uniform randomly selected initial conditions. 
The simulation results show that Algorithm~\ref{ADMMalgorithm} is about two orders of magnitude faster than the one in \cite{salehisadaghianiPavel2018Auto}. 
\vspace{-0.25cm}

\section{Conclusions}
\vspace{-0.3cm}We designed a distributed NE seeking algorithm in general games 
by using an inexact-ADMM approach. Each player maintains action estimates for the others, exchanged with its neighbours over a connected, undirected communication graph. The game was reformulated as  a modified game 
with virtual consensus constraints for the action estimates, and solved within the framework of inexact-ADMM.  An inexact-ADMM algorithm was designed  which was shown to convergence to the NE of the game {under strong monotonicity of the pseudo-gradient and Lipschitz continuity of the extended pseudo-gradient}.  The convergence rate of the algorithm was compared with that of an existing gossip-based NE seeking algorithm.
{To deal with the merely monotone case (thus with possibly multiple NEs) we leave for future work; this would require more involved computations even on the consensus subspace. } As other directions  for future work we can mention: (1) extension to directed communication graphs, (2) adaptation to games with networked structure (where computation and memory could be reduced if each player maintains only its  own required estimates), (3) designing incentives to ensure that players
truthfully share their states and auxiliary variables.  
\vspace{-0.255cm}


\appendix 



\vspace{-0.3cm}
\section{Derivation of  ADMM Algorithm \ref{ADMMalgorithm}}
\vspace{-0.3cm}
In the following we show how  Algorithm \ref{ADMMalgorithm} can be derived based on an ADMM approach adapted for a game setup. 
Using slack variables $t^{ij}$, $t^{ji}$ associated with each edge to separate the constraints and $\mathcal{I}_{\Omega_i}(x_i^i)$ for the feasibility constraint $x_i^i\in\Omega_i$, we write \eqref{mini_1} as
 \vspace{-0.25cm}
\begin{eqnarray}
\label{mini_2}
&&\hspace{-0.09cm}\begin{cases}
\begin{aligned}
& \underset{x_i^i\in \mathbb{R}}{\text{minimize}}
& & J_i(x_i^i,x_{-i}^i)+\mathcal{I}_{\Omega_i}(x_i^i), \qquad \qquad  \forall i\in V \\ 
& \text{subject to}
&&  \hspace{-0.2cm} x^i=t^{ij};  \, \forall  j \, \text{ s.t. }  (i, j) \in E  \\  
&& &  \hspace{-0.2cm}  x^i =t^{ji};  \, \forall j \, \text{ s.t. }  (j,i) \in E 	
\end{aligned}
\end{cases}
\end{eqnarray} 
and similarly for player $j$.
We cannot  directly use a standard ADMM algorithm because each minimization in \eqref{mini_2} is not over the whole $x^i$ as in DOPT, cf. \eqref{miniBis}, but rather only over part of it (action $x^i_i$), and inherent coupling to the others' decisions still exists indirectly, via the estimates $x_{-i}^i$. 
\emph{Note that in  \eqref{mini_2},  since minimization of $J_i$ is only over the action component $x^i_i$,  for the other (estimate) components $x_{-i}^i$ we can consider a zero objective function}. We will develop an ADMM-type algorithm to solve \eqref{mini_2} in a distributed manner.   {Let $\{u^{ij},v^{ji}\} \in \mathbb{R}^N$  be the Lagrange multipliers of player $i$, associated with the  constraints in \eqref{mini_2} (source, destination), and similarly,  $\{u^{ji},v^{ij}\} \in \mathbb{R}^N$,  for player $j$, respectively}. 
Let   $N^{out}_i =\{ j | \,(i,j)\in E\}$ and $N^{in}_i=\{ j | \, (j,i)\in E\}$, and since $G_c$ is undirected, $N^{in}_i=N^{out}_i= N_i$. 
We show next how the algorithm can be derived  and simplified so that each player maintains only its dual variables $w^i:=\sum_{j\in N_i}u^{ij}+v^{ji} $, not all its Lagrangian multipliers $\{u^{ij},v^{ji}\}$. The auxiliary variables $t^{ij}$, $t^{ji}$ are obtained  based on minimizing the overall augmented Lagrangian, 
{\begin{eqnarray}\label{augmented_Lag_all} 
&&\hspace{-0.15cm}\mathcal{L}^{c}\big(\!\{x^i\},\! \{t^{ij},t^{ji}\}\}; \!\{u^{ij},v^{ji}\};\! \{u^{ji},v^{ij}\}\!\big) \!:= \! \sum_{i\in V}\!J_i(x_i^i,x_{-i}^i) 
 \nonumber\\
&&\hspace{-0.15cm}+ \sum_{i\in V}\sum_{j \in N_i}  {u^{ij}}^T(x^i-t^{ij}) 
+  \sum_{i\in V} \sum_{j \in N_i}  {v^{ji}}^T(x^i-t^{ji}) \\  
&& \hspace{-0.15cm}+\frac{c}{2} \sum_{i\in V} \sum_{j \in N_i}   \| x^i-t^{ij}\|^2  
+\frac{c}{2} \sum_{i\in V}\sum_{j \in N_i}  \| x^i-t^{ji}\|^2  \nonumber  
\end{eqnarray}}
Based on \eqref{augmented_Lag_all}, $t^{ij}$, $t^{ji}$  are updated  as: 
\begin{eqnarray*}
&&\hspace{-0.15cm} t^{ij}(k)\!=\!\text{arg }\min_{t^{ij}} \!\mathcal{L}^{c}\!\big(\!\{x^i(k\!-\!1)\},\! \{t^{ij}\}; \{\!u^{ij}\!(k\!-\!1)\!,\!v^{ij}\!(k\!-\!1\!)\!\}\!\big)\\
&&\hspace{-0.15cm} t^{ji}(k)\!=\!\text{arg }\min_{t^{ji}} \!\mathcal{L}^{c}\!\big(\!\{x^j(k\!-\!1)\},\! \{t^{ji}\}; \{\!u^{ji}\!(k\!-\!1)\!,\!v^{ji}\!(k\!-\!1\!)\!\}\!\big)
\end{eqnarray*}
or, since $\mathcal{L}^c$ is quadratic in $t^{ij}$, $t^{ji}$, in closed-form as  \vspace{-0.4cm}
\begin{eqnarray}\label{tij_asli}
&& \hspace{-0.15cm} t^{ij}(k)\!=\!\frac{\!u^{ij}(k\!-\!1)\!+\!v^{ij}(k\!-\!1)\!}{2c}\!+\!\frac{\!x^i(k\!-\!1)\!+\!x^j(k\!-\!1)\! }{2}\\
&& \hspace{-0.15cm} t^{ji}(k)\!=\!\frac{\!u^{ji}(k\!-\!1)\!+\!v^{ji}(k\!-\!1)\!}{2c}\!\!+\!\frac{\!x^j(k\!-\!1)\!+\!x^i(k\!-\!1)\!}{2}\! \nonumber
\end{eqnarray}
For each player $i\in V$ we consider an augmented Lagrangian associated to \eqref{mini_2}, where $c>0$, \vspace{-0.3cm} 
\begin{eqnarray}\label{augmented_Lag}
&&\hspace{-0.15cm}\mathcal{L}_i^{c}\big(x^i, t^{ij}, t^{ji};{ \{u^{ij},v^{ji}\}}\big):= J_i(x_i^i,x_{-i}^i)+\mathcal{I}_{\Omega_i}(x_i^i) \nonumber\\
&&\hspace{-0.15cm}+\sum_{j \in N_i}  {u^{ij}}^T(x^i-t^{ij})  
+  \sum_{j \in N_i} {v^{ji}}^T(x^i-t^{ji}) \\ 
&& \hspace{-0.15cm}+\frac{c}{2}\sum_{j \in N_i}  \| x^i-t^{ij}\|^2 
+\frac{c}{2}\sum_{j \in N_i}  \| x^i-t^{ji}\|^2  \nonumber 
\end{eqnarray}
For player $i$  {its dual variables $\{u^{ij}, v^{ji}\}$} are updated by dual ascent  with $c$ as the step-size, based on \eqref{augmented_Lag},  \vspace{-0.3cm} 
\begin{eqnarray}
&&\hspace{0.3cm}u^{ij}(k)=u^{ij}(k-1)+c\big(x^i(k-1)-t^{ij}(k)\big),\label{uijbefore}\\
&&\hspace{0.3cm}v^{ji}(k)=v^{ji}(k-1)+c\big(x^i(k-1)-t^{ji}(k)\big). \nonumber 
\end{eqnarray}
{Similarly, for player $j$'s, based on $\mathcal{L}_j^{c}$, its  $\{u^{ji}, v^{ij}\}$} are,\vspace{-0.3cm}
\begin{eqnarray}
&&\hspace{0.3cm}u^{ji}(k)=u^{ji}(k-1)+c\big(x^j(k-1)-t^{ji}(k)\big), \label{vijbefore}\\
&&\hspace{0.3cm}v^{ij}(k)=v^{ij}(k-1)+c\big(x^j(k-1)-t^{ij}(k)\big). \nonumber 
\end{eqnarray}
Considering  $u^{ij}(0)=v^{ij}(0)=u^{ji}(0)=v^{ji}(0)=\textbf{0}_N$ and using  \eqref{tij_asli} in  \eqref{uijbefore},\eqref{vijbefore} yields 
\vspace{-0.3cm}
\begin{eqnarray*}\label{u+v=0}
&&\hspace{0.3cm}u^{ij}(k)+v^{ij}(k)=\textbf{0}_N\quad \forall \,k\geq0\\ 
&&\hspace{0.3cm}u^{ji}(k)+v^{ji}(k)=\textbf{0}_N\quad \forall \,k\geq0. \nonumber  
\end{eqnarray*}
so that  \eqref{tij_asli} become, \vspace{-0.3cm} 
\begin{equation}
t^{ij}(k)=\frac{x^i(k-1)+x^j(k-1)}{2} =  t^{ji}(k).\label{t_ij_simp}
\end{equation}
Using  \eqref{t_ij_simp} to eliminate $t^{ij}(k)$, $t^{ji}(k)$ from \eqref{uijbefore}, yields \vspace{-0.3cm} 
\begin{eqnarray*}
&&\hspace{0.3cm}u^{ij}(k)=u^{ij}(k-1)+\frac{c}{2}\big(x^i(k-1)-x^j(k-1)\big) \label{u_update_ADMM}\\
&&\hspace{0.3cm}v^{ji}(k)=v^{ji}(k-1)+\frac{c}{2}\big(x^i(k-1)-x^j(k-1)\big).\nonumber 
\end{eqnarray*}
Then for  $w^i(k):=\sum_{j\in N_i}u^{ij}(k)+v^{ji} (k) $, this yields  \vspace{-0.3cm} 
\begin{equation}\label{w^i(k)}
w^i(k)=w^i(k-1)+c\sum_{j\in N_i}(x^i(k-1)-x^j(k-1)) 
\end{equation}
where $c>0$, as in Step 5 of Algorithm \ref{ADMMalgorithm}. Thus, { the auxiliary variables $t^{ij}$, $t^{ji}$ can be eliminated using \eqref{t_ij_simp} and,  
based on \eqref{augmented_Lag},  for each player $i\in V$  only the dual variables $w^i(k)$ \eqref{w^i(k)} are needed to update his \emph{action} $x_i^i(k)$ and \emph{estimates} $x_{-i}^i(k)$, using his and his neighbours' previous estimates, $x^i(k-1)$, $x^j(k-1)$, $ j\in N_i$. Specifically, }   
 minimizing \eqref{augmented_Lag}  w.r.t. $x^i_i$ and  with \eqref{t_ij_simp}, for each $i \in V$ 
the local update for $x_i^i(k)=x_i(k)$ (decision)  is \vspace{-0.25cm}
%
\begin{eqnarray}\label{x_i^i_ADMM_Pre}
&&\hspace{-0.15cm}x_i^i(k)=\text{arg }\min_{x_i^i\in\mathbb{R}}\Big\{J_i(x_i^i,{x_{-i}^i(k-1)})+\mathcal{I}_{\Omega_i}(x_i^i)\\		
&&\hspace{0.5cm}+w_i^i(k)^Tx_i^i+c\sum_{j\in N_i}\Big\|x_i^i-\frac{x_i^i(k-1)+x_i^j(k-1)}{2}\Big\|^2\Big\}.\nonumber
\end{eqnarray}
{Note that in DOP, where the cost $f_i(x_i)$ is separable in $x_i$ (cf. \eqref{miniBis}), in the update of $x_i$,  coupling appears only linearly due to the quadratic terms, e.g., \cite{WeiOzdaglar2013arX}, \cite{HongLuoRaza2016}. In contrast, in \eqref{x_i^i_ADMM_Pre} the $x_i^i$ problem is coupled nonlinearly to $x_{-i}^i(k-1)$ (due to the inherent coupling in $J_i$) and linearly to $x_i^j(k-1)$.} 
As in an \emph{inexact ADMM},   instead of solving exactly an optimization sub-problem, each player uses a linear proximal approximation to reduce the complexity of each sub-problem. We approximate $J_i(x_i^i,x_{-i}^i(k-1))$ in \eqref{x_i^i_ADMM_Pre} by a proximal first-order approximation  around $x^i_i(k-1)$, but {inspired by  \cite{HongLuoRaza2016}, with an extra quadratic penalty}, which yields \vspace{-0.3cm}
\begin{eqnarray}\label{x_i^i_ADMM_Pre1}
&&\!x_i^i(k\!)\!=\!\text{arg }\!\!\min_{x_i^i\in\mathbb{R}}\!\! \left \{\nabla_iJ_i(x^i(k\!-\!1))^T\!(x_i^i\!  \right .
  -  \!x_{i}^i(k\!-\!1\!)) \nonumber\\
 &&\hspace{0.8cm}\!+\!\frac{\beta_i}{2}\|x_i^i\! -\!x_{i}^i(k\!-\!1)\|^2\! +\!\mathcal{I}_{\Omega_i}(x_i^i)\\
&&\hspace{0.8cm}\! \left  . +w_i^i(k)^Tx_i^i\!+\!{\bar{c}} \sum_{j\in N_i}\Big\|x_i^i\!-\!\frac{x_i^i(k-1)\!+\!x_i^j(k\!-\!1)}{2}\Big\|^2\right\},\nonumber 
\end{eqnarray}
where   $\beta_i\!>\!0$, {$\bar{c} =c+ c_0$, $c_0\geq0$}. This is equivalent to \vspace{-0.3cm}
\begin{eqnarray}\label{x_i^i_ADMM_Pre2}
&&\hspace{-0.15cm}\!x_i^i(k)\!=\!\text{arg }\!\!\min_{x_i^i\in\mathbb{R}}\!\Big\{\mathcal{I}_{\Omega_i}(x_i^i)\!+\!\frac{\alpha_i}{2}\Big\|\!x_i^i\!-\!\alpha_i^{-1} \!\Big [
-\!\nabla_iJ_i(x^i(k\!-\!1)\!) \! \nonumber \\
&&\hspace{0.4cm}
+ \beta_ix_i^i(k\!-\!1\!) \!-\!w_i^i(k)\!+\! {\bar{c}}  \sum_{j\in N_i}\!\big(\!x_i^i(k\!-\!1\!)\!+\!x_i^j(k\!-\!1\!)\big)\! \Big ] \!\Big\|^2\Big\}, \nonumber 
\end{eqnarray}
where $\alpha_i =\beta_i + 2 {\bar{c} } |N_i|$. 
Using 
 \eqref{prox_def} and $\text{prox}_{\mathcal{I}_{\Omega_i}}^{\alpha_i}\{ \cdot \} = T_{\Omega_i} \{ \cdot \}  $, 
 each player $i$ updates its \emph{action} $x^i_i(k) = x_i(k)$ 
as in Step 6 of Algorithm \ref{ADMMalgorithm}. 
Finally, 
 each player $i$ updates its \emph{estimates} of the other players' actions $x_{-i}^i(k)$. This update is similar to \eqref{x_i^i_ADMM_Pre} 
  but ignoring the $J_i(x_i^i,x_{-i}^i)$ term (which is minimized only w.r.t. $x^i_i$, see \eqref{mini_2}), i.e.,   \vspace{-0.4cm}
\begin{eqnarray*}\label{x_i^i_ADMM_Pre-i}
&&\hspace{-0.15cm}x_{-i}^i(k)=\text{arg }\min_{x_{-i}^i\in\mathbb{R}^{N-1}}\Big\{w_{-i}^i(k)^Tx_{-i}^i \\
&&\hspace{0.5cm} +\bar{c}\sum_{j\in N_i}\Big\|x_{-i}^i-\frac{x_{-i}^i(k-1)+x_{-i}^j(k-1)}{2}\Big\|^2\Big\}.
\end{eqnarray*}  
Following the same proximal approximation, 
 the update for   $x^i_{-i}(k)$ is obtained  
as in Step 7 of Algorithm \ref{ADMMalgorithm}. 

\vspace{-0.25cm}

\section{Proof of Proposition \ref{propAlg1_vector}}
\vspace{-0.25cm}
The stacked vector update for $\mathbf{w}$,  \eqref{wstep_Algvector},  follows from Step~5 of Algorithm~\ref{ADMMalgorithm}. 
{From Step~6 in Algorithm~\ref{ADMMalgorithm},    $\alpha_i x_i^i(k) = T_{\Omega_i} \{ y \} = \text{prox}_{\mathcal{I}_{\Omega_i}}\{y \}$, where $y$ is  the term in the square-bracket on the right-hand side of Step~6.  
From Proposition 16.34, \cite{Combettes},  
$p  =\text{prox}_{\mathcal{I}_{\Omega_i}}\{ y \}$ $\Leftrightarrow$ $y-p  \in  \partial \mathcal{I}_{\Omega_i}(p )$. Thus,  $\alpha_i x_i^i(k)  =\text{prox}_{\mathcal{I}_{\Omega_i}}\{ y \}$ $\Leftrightarrow$ $y-\alpha_i x_i^i(k)  \in  \partial \mathcal{I}_{\Omega_i}(\alpha_i x_i^i(k) )$ $\Leftrightarrow$ $y-\alpha_i x_i^i(k)  \in  \partial \mathcal{I}_{\Omega_i}(x_i^i(k) )$ since $\alpha_i >0$. Expanding  $y$ as on the right-hand side of Step~6 and using  $\alpha_i = \beta_i + 2 {\bar{c}} | N_i |>0$  this is equivalent to \vspace{-0.3cm} 
\begin{eqnarray*}\label{step6_alt}&&\hspace{-0.15cm} \mathbf{0} { \in} \,  \nabla_iJ_i(x_i^i(k-1),x^i_{-i}(k-1)) \!+\!\partial_i\mathcal{I}_{\Omega_i}(x_i^i(k)) \!+\!w_i^i(k) \!+\! \nonumber \\ &&\hspace{-0.15cm}  (\beta_i\!+ \!2 {\bar{c}} |N_i|)x_i^i(k\!)  \!-\! \beta_ix_i^i(k\!-\!1\!) \!-\! \!{\bar{c}}\!\sum_{j\in N_i}\!\!(\!x_i^i(k\!-\!1\!)\!+\!x_i^j(k\!-\!1\!)) 
 \end{eqnarray*} \vspace{-0.2cm}i.e., } \vspace{-0.2cm}
 \begin{eqnarray*}\label{step6_alt}
&&\hspace{-0.15cm} \mathbf{0} { \in} \,  \nabla_iJ_i(x_i^i(k-1),x^i_{-i}(k-1))\! +\!\partial_i\mathcal{I}_{\Omega_i}(x_i^i(k)) \!+\!w_i^i(k)+  \nonumber \\
&&\hspace{-0.15cm} (\beta_i\!+\! 2 {\bar{c}} |N_i|\!)(x_i^i(k)  \!-\!x_i^i(k\!-\!1)\!) \!+\! { \bar{c} }\!\sum_{j\in N_i}\!\!(\!x_i^i(k\!-\!1\!)\!-\!x_i^j(k\!-\!1\!)).
\end{eqnarray*}	
From Step~7 in Algorithm~\ref{ADMMalgorithm}, one can obtain for  $x_{-i}^i$, \vspace{-0.3cm} 
\begin{eqnarray*}\label{step7_alt}
&&\hspace{-0.15cm} \mathbf{0}_{N-1} \, { =} \, (\beta_i +2 {\bar{c}} |N_i| )(x_{-i}^i(k) -  x_{-i}^i(k-1))   +w_{-i}^i(k) \\
&&\hspace{0.8cm}  + { \bar{c}} \! \sum_{j\in N_i} (x_{-i}^i(k\!-\!1) \!-\! x_{-i}^j(k\!-\!1)). \nonumber  
\end{eqnarray*}	
Combining  these two relations  into a single vector one for all components of $x^i = (x_i^i, x_{-i}^i)$, yields $ \forall i\in V $ \vspace{-0.3cm} 
\begin{eqnarray*}
&&\hspace{-0.15cm}  \mathbf{0}_N \, { \in} \,  \nabla_iJ_i(x^i(k-1)) \, e_i  + \partial_i\mathcal{I}_{\Omega_i}(x_i^i(k)) \, e_i \!+\!w^i(k) \!+ \nonumber \\
&&\hspace{-0.2cm} \!(\beta_i\!+\! 2 {\bar{c}} |N_i|\!)(x_i^i(k)  \!-\!x_i^i(k\!-\!1)\!) \!+\! { \bar{c} }\!\sum_{j\in N_i}\!\!(\!x_i^i(k\!-\!1\!)\!-\!x_i^j(k
\!-\!1\!)).
\end{eqnarray*}	
In stacked vector form for all $i \in V$  this yields \vspace{-0.3cm}
\begin{eqnarray*}\label{step67_altVec00}
&&\hspace{-0.15cm} \!\mathbf{0}_{N^2}\, { \!\in} \! \,  [ \nabla_iJ_i(x^i(k\!-\!1)\!)  e_i\!]_{i \in V} \! +\! [\partial_i\mathcal{I}_{\Omega_i}(x_i^i(k)) e_i \!]_{i \in V} \!+ \!\mathbf{w}(k)  \nonumber \\
&&\hspace{0.5cm}+((\mathbb{B} +2 {\bar{c}} \mathbb{D} )\otimes I_N) (\mathbf{x}(k) - \mathbf{x}(k-1))  \nonumber\\
&&\hspace{0.5cm}  + { \bar{c}} (L \otimes I_N) \mathbf{x}(k-1),  
\end{eqnarray*}	
where $\mathbf{x} = [x^i]_{i \in V}$, $\mathbf{w} = [w^i]_{i \in V}$,  $ \mathbb{B} = \text{diag}([\beta_i]_{i \in V})$, $\mathbb{D} = \text{diag}([| N_i  |]_{i \in V})$, 
which  can be written as  \eqref{xstep_step67_altVec0}. 
$\hfill\blacksquare$\vspace{-0.25cm}

{{\bf Proof of Lemma \ref{boldRF_plus_cL_strgmon_lemma}}} \vspace{-0.4cm}

{
We decompose  the augmented space $ \mathbb{R}^{N^2}$ into the  consensus subspace $Ker(\mathbf{L}) =\{ \mathbf{1}_N\otimes x  \, | x \in \ \mathbb{R}^{N} \}$ and its orthogonal complement  $Ker(\mathbf{L})^{\perp}$. 
Any $\mathbf{x} \in \mathbb{R}^{N^2}$  
can be decomposed as $\mathbf{x}= \mathbf{x}^{\|} + \mathbf{x}^{\perp}$, with $ \mathbf{x}^{\|} \in  Ker (\mathbf{L})$, $\mathbf{x}^{\perp} \in  Ker (\mathbf{L})^{\perp}$,  $ ( \mathbf{x}^{\|})^T\mathbf{x}^{\perp} =0$, 
Thus $\mathbf{x}^{\|} = \mathbf{1}_N \otimes x$, for some $x \in  \mathbb{R}^{N}`$, 
so that  $\mathbf{L} \mathbf{x}^{\|} =0$, and $min_{\mathbf{x}^{\perp}\in Ker(\mathbf{L})^{\perp}} (\mathbf{x}^{\perp})^T \mathbf{L} \mathbf{x}^{\perp} = \lambda_2(L) \| \mathbf{x}^{\perp} \|^2$, where $ \lambda_2(L) >0$.
Similarly, we can write any $\mathbf{y} \in Ker(\mathbf{L})$ as $\mathbf{y}= \mathbf{1}_N \otimes y$, for some $y \in  \mathbb{R}^{N}$  
and $\mathbf{L} \mathbf{y}=0$. 
Then, using  $\mathbf{x} = \mathbf{x}^{\|} + \mathbf{x}^{\perp}$, $\mathbf{F}(\mathbf{x}^{\|}) = F(x)$, $\mathbf{F}(\mathbf{y}) = F(y)$, $\mathbf{R}^T \mathbf{x}^{\|} = x$, $\mathbf{R}^T \mathbf{y} = y$ we can write \vspace{-0.25cm}
\begin{eqnarray}\label{eqboldRF_plus_cL_strgmon_1}
&&(\!\mathbf{x}\!-\! \mathbf{y})^T \! \big ( \mathbf{R} (\mathbf{F}(\mathbf{x})\!-\!\mathbf{F}(\mathbf{y}))\!+\!c_0 \mathbf{L}(\mathbf{x} \!-\! \mathbf{y}) \big )\!  \nonumber \\
&&\hspace{0.5cm} = (\!\mathbf{x}^{\|}\!-\! \mathbf{y})^T\mathbf{R} \big (\mathbf{F}(\mathbf{x})\!-\!\mathbf{F}(\mathbf{x}^{\|})\! +\! \mathbf{F}(\mathbf{x}^{\|})\!  -\!\mathbf{F}(\mathbf{y}) \big)\! \nonumber \\
&&\hspace{0.8cm} + (\mathbf{x}^{\perp})^T\mathbf{R} \big ( \mathbf{F}(\mathbf{x})\!-\!\mathbf{F}(\mathbf{x}^{\|})\! +\! \mathbf{F}(\mathbf{x}^{\|})\!  -\!\mathbf{F}(\mathbf{y}) \big )\! \nonumber \\
&&\hspace{0.8cm} + c_0  (\mathbf{x}^{\|}\!+\!\mathbf{x}^{\perp}-\! \mathbf{y})^T\mathbf{L}(\mathbf{x}^{\|} + \mathbf{x}^{\perp} \!-\! \mathbf{y}) \big ) \! \nonumber \\
&&\hspace{0.5cm} = (x-y)^T \!\big(\mathbf{F}(\mathbf{x})\!-\!\mathbf{F}(\mathbf{x}^{\|}) \big )\! +  (x-y)^T \! \big ( F(x)\!  -F(y) \big )\! \nonumber \\
&&\hspace{0.8cm} + (\!\mathbf{x}^{\perp})^T \mathbf{R} \big (\mathbf{F}(\mathbf{x})\!-\!\mathbf{F}(\mathbf{x}^{\|}) \big )\! + \! (\!\mathbf{x}^{\perp})^T \mathbf{R}  \big (F(x)\!  -\! F(y) \big )\! \nonumber \\
&&\hspace{0.8cm} +c_0 (\mathbf{x}^{\perp})^T \mathbf{L} \mathbf{x}^{\perp} 
\end{eqnarray}
Using strong monotonicity of $F$ (by Assumption \ref{strgmon_Fassump}) for the second term, and properties of $ \mathbf{L} $ on $Ker(\mathbf{L})^{\perp} $ for the fifth term  on the right-hand side we can write,  \vspace{-0.25cm}
\begin{eqnarray}\label{eqboldRF_plus_cL_strgmon_2}
&&(\!\mathbf{x}\!-\! \mathbf{y})^T \! \big ( \mathbf{R} (\mathbf{F}(\mathbf{x})\!-\!\mathbf{F}(\mathbf{y}))\!+\!c_0 \mathbf{L}(\mathbf{x} \!-\! \mathbf{y}) \big )\!  \nonumber \\
&&\hspace{0.5cm} \geq (x-y)^T \!\big(\mathbf{F}(\mathbf{x})\!-\!\mathbf{F}(\mathbf{x}^{\|}) \big )\! + \mu \| x-y \|^2 \nonumber \\
&&\hspace{0.8cm} + (\!\mathbf{x}^{\perp})^T \mathbf{R} \big (\mathbf{F}(\mathbf{x})\!-\!\mathbf{F}(\mathbf{x}^{\|}) \big )\! + \! (\!\mathbf{x}^{\perp})^T \mathbf{R}  \big (F(x)\!  -\! F(y) \big )\! \nonumber \\
&&\hspace{0.8cm} + c_0  \lambda_2 (L) \| \mathbf{x}^{\perp} \|^2
\end{eqnarray}
We deal with the cross-terms by using $a^T b \geq - \|a\| \|b \|$, for any $a$, $b$, and $\| \mathbf{R} b \| \leq \| \mathbf{R}\|_2 \|b\| $, so that  
\vspace{-0.25cm}
\begin{eqnarray}\label{eqboldRF_plus_cL_strgmon_3}
&&(\!\mathbf{x}\!-\! \mathbf{y})^T \! \big ( \mathbf{R} (\mathbf{F}(\mathbf{x})\!-\!\mathbf{F}(\mathbf{y}))\!+\!c_0 \mathbf{L}(\mathbf{x} \!-\! \mathbf{y}) \big )\! \geq	 \nonumber \\
&&\hspace{0.3cm}  - \|x-y\| \, \| \mathbf{F}(\mathbf{x})\!-\!\mathbf{F}(\mathbf{x}^{\|}) \| \! + \mu \| x-y \|^2 \nonumber \\
&&\hspace{0.3cm} - \|\!\mathbf{x}^{\perp} \| \|\mathbf{R}\|_2 \| \mathbf{F}(\mathbf{x})\!-\!\mathbf{F}(\mathbf{x}^{\|}) \| \! - \! \| \!\mathbf{x}^{\perp}\| \|\mathbf{R}\|_2 
\| F(x)\!  -\! F(y) \|\! \nonumber \\
&&\hspace{0.3cm} + c_0  \lambda_2 (L) \| \mathbf{x}^{\perp} \|^2
\end{eqnarray}
Using  $\| \mathbf{R}\|_2 = 1$ and  Lipschitz properties of $F$ and $\mathbf{F}$,  \vspace{-0.3cm} 
\begin{eqnarray}\label{eqboldRF_plus_cL_strgmon_4}
&&(\!\mathbf{x}\!-\! \mathbf{y})^T \! \big ( \mathbf{R} (\mathbf{F}(\mathbf{x})\!-\!\mathbf{F}(\mathbf{y}))\!+\!c_0 \mathbf{L}(\mathbf{x} \!-\! \mathbf{y}) \big )\!  \geq \nonumber \\
&&\hspace{0.3cm}  - \theta \|x-y\| \, \| \mathbf{x}^{\perp} \| \! + \mu \| x-y \|^2 - \theta \| \mathbf{x}^{\perp} \|^2 \!  - \theta_0 \! \| \mathbf{x}^{\perp}\| 
\| x\!  -\! y \|\!  \nonumber \\
&&\hspace{0.3cm} 
 + c_0  \lambda_2 (L) \| \mathbf{x}^{\perp} \|^2
\end{eqnarray}
Using $ \|x-y\| = \frac{1}{\sqrt{N}}  \| \mathbf{x}^{\|} \!-\! \mathbf{y}\|$ and $\Psi $, \eqref{eqboldRF_plus_cL_strgmon_5}, we can write 
\vspace{-0.25cm}
\begin{eqnarray*}
&&\hspace{-0.15cm}(\!\mathbf{x}\!-\! \mathbf{y})^T \!\! \big ( \!\mathbf{R} (\mathbf{F}(\mathbf{x}\!)\!-\!\mathbf{F}(\mathbf{y})\!)\!+\!c_0 \mathbf{L}(\mathbf{x} \!-\! \mathbf{y}\!)\! \big )\!  \geq \!
\!\left [ \! \begin{array}{c}
\! \| \mathbf{x} ^{\|}\!-\! \mathbf{y} \| \!\\
  \| \mathbf{x}^{\perp} \|
 \end{array} 
\!  \right ]^{\!T}
\! \!\! \!\! \Psi  \! \left [ 
\! \begin{array}{c}
 \! \| \mathbf{x} ^{\|}\!-\! \mathbf{y} \|\! \\
  \| \mathbf{x}^{\perp} \|
 \end{array} 
 \! \right ]
\end{eqnarray*}
Thus,  \vspace{-0.25cm}
$$
(\!\mathbf{x}\!-\! \mathbf{y})^T \! \big ( \mathbf{R} (\mathbf{F}(\mathbf{x}\!)\!-\!\mathbf{F}(\mathbf{y})\!)\!+\!c_0 \mathbf{L}(\mathbf{x} \!-\! \mathbf{y}\!) \big )\!  \geq \! \bar{\mu} 
\left (
 \| \mathbf{x} ^{\|}\!-\! \mathbf{y} \|^2 + 
  \| \mathbf{x}^{\perp} \|^2 \right), 
$$where $\bar{\mu}:=\lambda_{\min}(\Psi)$. Using $ \| \mathbf{x} ^{\|}\!-\! \mathbf{y} \|^2 +   \| \mathbf{x}^{\perp} \|^2  =  \| \mathbf{x} \!-\! \mathbf{y} \|^2$, this can be written as \eqref{eqboldRF_plus_cL_strgmon}. For $c_0$ as in the statement, $\Psi \succ 0$, hence $\bar{\mu}>0$. 
 $\hfill\blacksquare$ }

\bibliography{ref}

\end{document}